\documentclass[11pt,reqno]{amsart}

\usepackage[T1]{fontenc}
\usepackage{lmodern}
\usepackage{microtype}
\usepackage[a4paper,hmargin=25mm,vmargin=22mm,headheight=14pt,headsep=12pt]{geometry}
\usepackage{amsmath,amssymb,amsthm,mathtools}
\usepackage{enumitem}
\usepackage{hyperref}
\usepackage{aliascnt}
\usepackage[nameinlink,noabbrev,capitalise]{cleveref}
\AtBeginDocument{%
  \setlength{\abovedisplayskip}{8pt plus 2pt minus 2pt}%
  \setlength{\belowdisplayskip}{8pt plus 2pt minus 2pt}%
  \setlength{\abovedisplayshortskip}{3pt plus 2pt}%
  \setlength{\belowdisplayshortskip}{5pt plus 2pt minus 2pt}%
}
\hypersetup{
  hidelinks,
  pdftitle={Extensions at Nonzero p-Characters for Cartan-Type Lie Algebras},
  pdfauthor={Chao Ma},
  pdfcreator={Chao Ma},
  pdfsubject={Peirce descent for nonzero p-characters and the Yoneda algebra of the height-two Witt module in characteristic seven},
  pdfkeywords={Peirce decomposition, trace reconstruction, reduced enveloping algebra, p-character, Morita equivalence, matrix corner, restricted Lie algebra, Cartan type, radical, primitive idempotent, modular representation theory, Yoneda algebra, Tate resolution, Witt algebra, self-extension}
}

\setlist[enumerate]{leftmargin=2.2em,itemsep=0.15em,topsep=0.25em}
\setlist[itemize]{leftmargin=2.0em,itemsep=0.15em,topsep=0.25em}

\newtheorem{theorem}{Theorem}[section]
\newaliascnt{proposition}{theorem}
\newtheorem{proposition}[proposition]{Proposition}
\aliascntresetthe{proposition}
\newaliascnt{lemma}{theorem}
\newtheorem{lemma}[lemma]{Lemma}
\aliascntresetthe{lemma}
\newaliascnt{corollary}{theorem}
\newtheorem{corollary}[corollary]{Corollary}
\aliascntresetthe{corollary}
\theoremstyle{definition}
\newaliascnt{definition}{theorem}
\newtheorem{definition}[definition]{Definition}
\aliascntresetthe{definition}
\theoremstyle{remark}
\newaliascnt{remark}{theorem}
\newtheorem{remark}[remark]{Remark}
\aliascntresetthe{remark}
\theoremstyle{definition}
\newaliascnt{example}{theorem}
\newtheorem{example}[example]{Example}
\aliascntresetthe{example}

\newcommand{\kk}{\mathfrak{k}}
\newcommand{\Achi}{A_{\chi}}
\newcommand{\Mat}{\operatorname{Mat}}
\newcommand{\PT}{\mathrm{PT}}
\newcommand{\Jac}{\operatorname{Jac}}
\newcommand{\Soc}{\operatorname{Soc}}
\newcommand{\Tr}{\operatorname{Tr}}
\newcommand{\Hom}{\operatorname{Hom}}
\newcommand{\Ext}{\operatorname{Ext}}
\newcommand{\rad}{\operatorname{rad}}
\newcommand{\ad}{\operatorname{ad}}
\newcommand{\Prim}{\operatorname{Prim}}
\newcommand{\spn}{\operatorname{span}}
\newcommand{\Fp}{\mathbf F_p}

\newcommand{\vp}{v_p}
\newcommand{\nuP}{\nu_p}

\title[Extensions at nonzero $p$-characters]{Extensions at Nonzero $p$-Characters for Cartan-Type Lie Algebras}
\author{Chao Ma}
\thanks{Independent Researcher, London, United Kingdom. Email:
\texttt{raymond.ma@me.com}.}
\date{}

\begin{document}

\begin{abstract}
For restricted simple Cartan-type Lie algebras over an algebraically
closed field of characteristic
$p\ge5$, the reduced enveloping algebra at every nonzero
$p$-character is a matrix algebra over a terminal symmetric corner,
and this corner determines the extensions.
Central-socle traces recover its radical and its simple factors, and
the first radical layer determines the degree-one extensions.
For the Witt algebra $W(1;\mathbf1)$, at even height $2\le h\le p-3$, the
unique simple module has self-extension dimension
$\min\{p-h-2,h+2\}$.
At height two in characteristic seven, we present the Yoneda algebra
as an eight-dimensional symmetric algebra tensored with a polynomial
algebra on three central degree-two classes. An equivariant Tate
resolution and cocycle identities prove the presentation in all degrees.
\end{abstract}

\maketitle

\section{Introduction}

Let $k$ be an algebraically closed field of characteristic $p\ge5$.
For a restricted simple Cartan-type Lie algebra $\mathfrak g$ and
$\chi\in\mathfrak g^*$, write
\[
 u_\chi(\mathfrak g)=U(\mathfrak g)/
 \bigl(x^p-x^{[p]}-\chi(x)^p1:x\in\mathfrak g\bigr).
\]
We construct matrix factors at every nonzero character in the
Witt--Jacobson, Special, Hamiltonian and Contact families. Successive
Peirce corners reduce the extension problem to a symmetric algebra
whose radical and simple factors are recovered by central-socle traces.
The first radical layer gives $\Ext^1$, and the corner multiplication
also determines higher extensions.

Benkart and Feldvoss ask for the extensions between simple modules for
restricted simple Cartan-type Lie algebras at an arbitrary character
\cite[Problem~4]{BenkartFeldvoss2015}.
We determine these groups at every nonzero character from the
multiplication in a terminal symmetric corner.

\begin{theorem}\label{thm:intro-main}
Let
\[
 \mathfrak g\in\{W(n),\ S(n)^{(1)}\ (n\ge3),\
 H(2m)^{(2)},\ K(2m+1;\mathbf1)^{(1)}\ (m\ge1)\}
\]
and let $0\ne\chi\in\mathfrak g^*$. Then
$u_\chi(\mathfrak g)\cong\Mat_{p^r}(C_\chi)$ with $r\ge1$, where
$C_\chi$ is a symmetric algebra containing no unital copy of
$\Mat_p(k)$. Central-socle traces on $C_\chi$ recover its radical $J$
and the simple factors $S_\alpha\cong\Mat_{n_\alpha}(k)$ of $C_\chi/J$,
with identities $z_\alpha$. For the simple $C_\chi$-modules
$L_\alpha,L_\beta$,
\[
 \dim\Ext^1_{C_\chi}(L_\alpha,L_\beta)
 =\dim\Hom_{S_\beta-S_\alpha}
 \bigl(L_\beta\otimes L_\alpha^*,\ z_\beta(J/J^2)z_\alpha\bigr).
\]
The relative bar resolution gives the higher extension groups, and all
these data transport to $u_\chi(\mathfrak g)$ under Morita equivalence.
\end{theorem}

This is proved as \cref{thm:global}. The construction begins with an embedded $\Mat_p(k)$ and
the coefficient isomorphism $A\cong\Mat_p(eAe)$.

For a graded algebra
$\mathfrak g=\bigoplus_{i=-q}^{r}\mathfrak g_i$, put
$\mathfrak g_{(j)}=\bigoplus_{i=j}^{r}\mathfrak g_i$.
Benkart and Feldvoss \cite[(2.4)]{BenkartFeldvoss2015} define the height by
\[
 h(\chi)=
 \begin{cases}
 \min\{-q\le j\le r:\chi(\mathfrak g_{(j)})=0\},
       &\chi(\mathfrak g_r)=0,\\
 r+1,&\chi(\mathfrak g_r)\ne0.
 \end{cases}
\]
For the Witt algebra $W(1;\mathbf1)$, at every even height $2\le h\le p-3$,
\[
 \dim\Ext^1(L,L)=\min\{p-h-2,h+2\}.
\]
At height $p-3$ the terminal corner is $k[t]/(t^p)$.
At height two in characteristic seven, the unique simple module has
dimension $49$, its basic corner has dimension $343$, and
\[
 \Ext^*(L,L)\cong B\otimes k[\theta_3,\theta_4,\theta_5],
 \qquad |\theta_i|=2,
\]
where $B$ is the eight-dimensional symmetric algebra presented in
\cref{thm:w11-p7-yoneda}. Its Hilbert series is $(1+t)^3$, and its
degree-one square-zero locus is a line.
Thus $\dim\Ext^n(L,L)=(n+1)(n+2)/2$ for every $n\ge0$.
The proof uses an equivariant Tate resolution, central degree-two
operators, and cocycle identities.

Boe, Nakano and Wiesner determined the extensions of simple
$W(1;1)$-modules for $h(\chi)\in\{-1,0,1\}$ and for
$\chi(\mathfrak g_r)\ne0$, and proved a lower bound for the
self-extension dimension of the unique simple module at the other
heights \cite[\S4]{BoeNakanoWiesner2009}.
Rian obtained the extensions for $h(\chi)\in\{-1,0,1\}$ and
$\chi(\mathfrak g_r)\ne0$ by different methods
\cite{RianThesis2011,Rian2011}.
\Cref{thm:bnw-main} gives the self-extension dimension at even height
$2\le h\le p-3$. \Cref{thm:w11-p7-yoneda} gives the Yoneda algebra at
$p=7$, $h=2$, the first even height at which $W(1;1)_{(h)}$ is
nonabelian.

Write $\kk=K(2m+1;\mathbf1)^{(1)}$, $m\ge1$, and
$\Achi=u_\chi(\kk)$. The Contact bracket supplies the following first
factor; \cref{sec:wsh} gives the other three families.

\begin{theorem}[The first matrix factor]\label{thm:first-matrix}
Let $0\ne\chi\in\kk^*$.  Then $u_\chi(\kk)$ contains a unital copy of
$\Mat_p(k)$.  If $c$ spans $K_{-2}$ and
$V=K_{-1}$, then
\[
\chi(c)\ne0\implies \Mat_{p^{m+1}}(k)\subseteq \Achi,
\]
with the same unit, and
\[
\chi(c)=0,\quad \chi|_V\ne0
\implies \Mat_{p^m}(k)\subseteq \Achi
\]
with the same unit.
\end{theorem}

A corner is called terminal when it contains no unital $\Mat_p(k)$.
For $C/J(C)=\prod_\lambda\Mat_{n_\lambda}(k)$, put
$\nuP(C)=\min_\lambda v_p(n_\lambda)$; a terminal algebra has
$\nuP(C)=0$ and is called $p$-primitive.

\begin{theorem}[$p$-Peirce descent]\label{thm:uniform-descent}
Let $C_0$ be a finite-dimensional symmetric $k$-algebra.  There is a finite sequence of symmetric algebras
\[
 C_0,C_1,\ldots,C_r
\]
such that
\[
 C_j\cong \Mat_p(C_{j+1})\quad(0\le j<r),
\]
and, for every $0\le j\le r$,
\[
 C_j\text{ contains a unital copy of }\Mat_p(k)
 \quad\Longleftrightarrow\quad
 \tau_j|_{Z(C_j)\cap\Soc(C_j)}=0,
\]
where $\tau_j$ is a symmetrizing functional on $C_j$.  Moreover
\[
 r=\nuP(C_0),\qquad
 C_0\cong\Mat_{p^r}(C_r),\qquad
 \nuP(C_r)=0.
\]
The resulting Morita equivalence preserves blocks, simple modules,
projective covers, and all extension groups.
\end{theorem}

The same central-socle pairing governs reconstruction of the endpoint.

\begin{theorem}[Trace reconstruction at a $p$-primitive endpoint]\label{thm:trace-reconstruction}
Let $C$ be a $p$-primitive endpoint algebra, let $\tau$ be a symmetrizing
functional, and put
\[
 R(C)=Z(C)\cap\Soc(C).
\]
There exists $s_0\in R(C)$ such that
\[
 B_{s_0}(x,y)=\tau(xys_0)
\]
has radical $\Jac(C)$.  Thus the form determines $S=C/\Jac(C)$.  The map
\[
 s\longmapsto \ell_s,\qquad
 \ell_s(\bar x)=\tau(xs),
\]
identifies $R(C)$ with the trace space of $S$ and, after fixing $s_0$, with
$Z(S)$.  For a generic $s_1\in R(C)$, the squarefree minimal polynomial of
the corresponding element of $Z(S)$ and Lagrange interpolation give all
primitive central idempotents
\[
 z_\alpha\in\Prim Z(S),\qquad
 S_\alpha=z_\alpha S\cong\Mat_{n_\alpha}(k).
\]
If $G_C$ is a finite generating set obtained from the iterated Peirce
coefficients, define
\[
 \rho_\alpha(g)=z_\alpha(g+\Jac(C)),\qquad g\in G_C.
\]
Finally, with $J=\Jac(C)$ and $M=J/J^2$,
\[
 M_{\beta\alpha}=z_\beta Mz_\alpha
\]
is a direct sum of copies of the unique simple
$(S_\beta,S_\alpha)$-bimodule $L_\beta\otimes L_\alpha^*$, and
\[
 \dim_k\Ext_C^1(L_\alpha,L_\beta)
 =\dim_k\Hom_{S_\beta-S_\alpha}
 (L_\beta\otimes L_\alpha^*,M_{\beta\alpha}).
\]
The normalized relative bar resolution computes the higher $\Ext$ groups.
\end{theorem}

Sections~\ref{sec:contact-setup}--\ref{sec:wsh} construct the matrix
factors. Sections~\ref{sec:peirce}--\ref{sec:main} prove the descent
and the trace reconstruction. Section~\ref{sec:harvest} gives the
basic algebra, and Section~\ref{sec:w11} treats the Witt algebra
$W(1;\mathbf1)$.

\section{The Contact algebra and matrix lemmas}\label{sec:contact-setup}

Assume $m\ge1$.  Let $\mathcal O$ be the height-one divided-power algebra in
the symplectic variables
\[
 u_1,\ldots,u_m,v_1,\ldots,v_m
\]
and the Contact variable $x_n$.  Put $V=\spn\{u_i,v_i\}_{i=1}^m$ and let
$E$ denote the Euler operator in the $V$-variables.  In the normalization
used below, the Contact bracket is
\begin{equation}\label{eq:contact-bracket}
 \{f,g\}_K
 =\{f,g\}_{\mathrm{sp}}
 +(2f-Ef)\partial_{x_n}g
 +(\partial_{x_n}f)(Eg-2g).
\end{equation}
Equivalently,
\begin{align*}
 \{f,g\}_K
 &=\Delta(f)D_ng-D_nf\,\Delta(g)
   +\sum_{j=1}^{2m}\sigma(j)D_jfD_{j'}g,\\
 \Delta(f)&=2f-\sum_{j=1}^{2m}x_jD_jf.
\end{align*}
Write $D_K(f)$ for the corresponding Contact derivation; see
\cite[Section~2]{Mygind2014}.  The derived simple restricted algebra has the
grading
\[
 \kk=\bigoplus_jK_j,
 \qquad K_{-2}=kc,\quad K_{-1}=V,
 \quad K_0\cong\mathfrak{sp}(V)\oplus kz,
\]
where $c=D_K(1)$ and $z$ is the grading element.  Normalize so that
\[
 [u_i,v_j]=\delta_{ij}c,
 \qquad [z,c]=-2c,
 \qquad [z,u_i]=-u_i,
 \qquad [z,v_i]=-v_i,
\]
and
\[
 c^{[p]}=u_i^{[p]}=v_i^{[p]}=0,
 \qquad z^{[p]}=z.
\]
If $\tau=(p-1,\ldots,p-1)$ denotes the full-top divided-power
exponent, then the divided-power Contact realization gives
\[
 K(2m+1;\mathbf1)^{(1)}
 =\spn\{D_K(x^{(\alpha)}):\alpha\ne\tau\}
 \quad\text{when }2m+4\equiv0\pmod p,
\]
and $K(2m+1;\mathbf1)^{(1)}=K(2m+1;\mathbf1)$ otherwise; see
\cite[Section~2]{Mygind2014}.

\begin{lemma}[Truncated Weyl algebra]\label{lem:weyl-matrix}
If $P,Q$ in a unital $k$-algebra satisfy
\[
 [P,Q]=1,\qquad P^p=Q^p=0,
\]
then the unital subalgebra they generate is $\Mat_p(k)$.
\end{lemma}

\begin{proof}
Let $Q$ act by multiplication by $t$ and $P$ by differentiation on
$k[t]/(t^p)$.  The relations hold and the resulting operators generate all
matrix units.  Conversely every word is a linear combination of
$Q^aP^b$, $0\le a,b<p$, so the universal algebra has dimension at most
$p^2$.  The displayed representation is onto a $p^2$-dimensional matrix
algebra and is therefore an isomorphism.  Any unital realization is the
image of a unital homomorphism from the simple algebra $\Mat_p(k)$, hence is
injective.
\end{proof}

\begin{lemma}[One-nilpotent Weyl factor]\label{lem:one-nilpotent}
If $P,Q$ in a unital $k$-algebra satisfy
\[
 [P,Q]=1,\qquad P^p=0,
\]
with no condition on $Q^p$, then the algebra contains a unital copy of
$\Mat_p(k)$.
\end{lemma}

\begin{proof}
Put $H=QP$.  The Weyl identity gives
\[
 H^p-H=Q^pP^p=0.
\]
For $a\in\Fp$ let
\[
 e_a=1-(H-a)^{p-1}.
\]
These are pairwise orthogonal idempotents with sum $1$.  If $e_0$ is the
zero-weight idempotent, set
\[
 E_{ij}=\frac1{j!}Q^ie_0P^j,
 \qquad 0\le i,j<p.
\]
The relations $[H,Q]=Q$ and $[H,P]=-P$ imply
\[
 e_0P^jQ^ke_0=\delta_{jk}j!e_0,
\]
so
\[
 E_{ij}E_{k\ell}=\delta_{jk}E_{i\ell}.
\]
Moreover
\[
 Q^ie_0P^i=i!e_i,
\]
whence $E_{ii}=e_i$ for every $0\le i<p$.  Since
$\sum_i e_i=1$, it follows that
\[
 \sum_iE_{ii}=1.
\]
Thus the $E_{ij}$ form a unital system of matrix units.
\end{proof}

\begin{lemma}[Finite Hasse slice]\label{lem:hasse-slice}
Let $A$ be an associative algebra of characteristic $p$, let
$\delta=\ad y$ with $\delta^p=0$ and $y^p=0$, and let $r\in A$ satisfy
$r^p=0$.  Put $u=\delta(r)$ and assume $u$ is a unit.  Define
\[
 w=\delta^{p-2}(r^{p-1}),
 \qquad v=\delta^{p-1}(r^{p-1}).
\]
If, after normalizing by the principal term $(p-1)!u^{p-1}$, every other
Hasse word lies in a nilpotent ideal of the finite algebra generated by the
Hasse string $r,\delta r,\ldots,\delta^{p-1}r$, then $v$ is a unit and
\[
 C=wv^{-1}
\]
satisfies $[y,C]=1$.  Consequently $A$ contains a unital
$\Mat_p(k)$.
\end{lemma}

\begin{proof}
The noncommutative Leibniz formula expands $v$ as
\[
 v=(p-1)!u^{p-1}+\eta,
\]
where $\eta$ is the sum of the nonprincipal Hasse words.  By hypothesis the
normalized error is nilpotent, so $v$ is a unit.  Since $\delta(v)=0$,
\[
 \delta(C)=\delta(w)v^{-1}=1.
\]
Thus $[y,C]=1$ and $y^p=0$.  Apply \cref{lem:one-nilpotent} with
$P=y$ and $Q=C$.
\end{proof}

\begin{lemma}[Invertible toral shift]\label{lem:toral-shift-strong}
Let $A$ be a unital, \emph{finite-dimensional} $k$-algebra and suppose
\[
 H^p-H=0,\qquad U\in A^\times,\qquad [H,U]=U.
\]
Then the unital subalgebra generated by $H,U$ contains a unital copy of
$\Mat_p(k)$.
\end{lemma}

\begin{proof}
The spectral idempotents $e_a=1-(H-a)^{p-1}$ satisfy
$e_ae_b=\delta_{ab}e_a$ and $\sum_ae_a=1$.
The relation $HU=U(H+1)$ gives
$Ue_aU^{-1}=e_{a+1}$. For $0\le i,j<p$ put
\[
 E_{ij}=U^ie_0U^{-j}.
\]
If $k\ne j$, then
$e_0U^{k-j}e_0=e_0e_{k-j}U^{k-j}=0$. Hence
\[
 E_{ij}E_{k\ell}=\delta_{jk}E_{i\ell},\qquad
 \sum_iE_{ii}=\sum_ie_i=1.
\]
The inverse $U^{-1}$ belongs to $k[U]$ by finite-dimensionality,
so all these matrix units lie in the stated subalgebra.
\end{proof}

\begin{corollary}[Scalar-centred version]\label{cor:toral-shift-strong}
Let $\mathfrak g$ be a finite-dimensional restricted Lie algebra.
If $h^{[p]}=h$ and
$[h,U]=U$ for a unit $U\in u_\chi(\mathfrak g)$, then a scalar translate
$H=h+\lambda1$ satisfies $H^p-H=0$ and
$u_\chi(\mathfrak g)$ contains a unital $\Mat_p(k)$.
\end{corollary}

\begin{proof}
Choose $\lambda$ with $\lambda^p-\lambda=-\chi(h)^p$.
In $u_\chi(\mathfrak g)$, $h^p=h+\chi(h)^p1$, so
$H^p-H=(h^p-h)+(\lambda^p-\lambda)1=0$, and $[H,U]=[h,U]=U$ since $\lambda$
is central.
\end{proof}

\begin{corollary}[Toral-shift matrix algebra]\label{lem:toral-shift}
If $h,u$ satisfy
\[
 [h,u]=u,\qquad h^p-h=0,\qquad u^p=1,
\]
then their unital subalgebra is $\Mat_p(k)$.
\end{corollary}

\begin{proof}
The relations span this subalgebra by $u^ah^b$, $0\le a,b<p$.
The matrix units of \cref{lem:toral-shift-strong} lie in it, so its
dimension is both at most and at least $p^2$.
\end{proof}

\begin{lemma}[Positive leading units]\label{lem:leading-unit}
Let $\mathfrak g$ be a restricted Lie algebra, graded by finitely many
integers, whose restricted $p$-map sends a homogeneous element of degree
$r>0$ to an element of degree $pr$ or to zero.  Let $r>0$ be maximal
with $\chi(\mathfrak g_r)\ne0$, and let $\Theta\in\mathfrak g_r$ satisfy
$\chi(\Theta)\ne0$.  Then $\Theta$ is a unit in $u_\chi(\mathfrak g)$.
\end{lemma}

\begin{proof}
Let $N_1=\Theta^{[p]}$.  Either $N_1=0$, or $N_1$ is homogeneous of
degree $pr>r$ and $\chi(N_1)=0$.  The same degree count applies to each
later restricted power, and finiteness of the grading forces
$\Theta^{[p^{q_0}]}=0$ for some $q_0$.  Whenever the character vanishes,
the reduced relation identifies the associative power with the restricted
power, so $N_1^{p^q}=\Theta^{[p^{q+1}]}$.  Thus $N_1$ is nilpotent, and
\[
 \Theta^p=\chi(\Theta)^p1+N_1
\]
is a nonzero scalar plus a nilpotent.
\end{proof}

\begin{lemma}[Nonzero leading weight]\label{lem:leading-weight}
Let $\mathfrak g$ and $\Theta$ be as in \cref{lem:leading-unit}, and
suppose $h\in\mathfrak g$ satisfies $h^{[p]}=h$ and
$[h,\Theta]=c\Theta$ for some $c\in\Fp^\times$.  Then
$u_\chi(\mathfrak g)$ contains a unital $\Mat_p(k)$.
\end{lemma}

\begin{proof}
\Cref{lem:leading-unit} gives that $\Theta$ is a unit.  Put
$h'=c^{-1}h$.  Then $[h',\Theta]=\Theta$.  Since $c^p=c$,
\[
 (h')^{[p]}=c^{-p}h^{[p]}=c^{-1}h=h'.
\]
Apply \cref{cor:toral-shift-strong}.
\end{proof}

\section{Negative-support matrix factors}\label{sec:negative}

\subsection{Central negative support}

\begin{proposition}\label{prop:central-negative}
If $a=\chi(c)\ne0$, then $\Achi$ contains a unital copy of
$\Mat_{p^{m+1}}(k)$.
\end{proposition}

\begin{proof}
The reduced relation gives $c^p=a^p1$, so $c$ is invertible.  Set
\[
 X_i=c^{-1}u_i,\qquad Y_i=v_i,
\]
and remove the scalar $p$-curvatures by setting
\[
 P_i=X_i-\frac{\chi(u_i)}a1,
 \qquad Q_i=Y_i-\chi(v_i)1.
\]
Then
\[
 [P_i,Q_j]=\delta_{ij},\qquad
 [P_i,P_j]=[Q_i,Q_j]=0,
 \qquad P_i^p=Q_i^p=0.
\]
By \cref{lem:weyl-matrix} the $m$ commuting Weyl pairs generate a unital
subalgebra
\[
 M_-\cong\Mat_{p^m}(k).
\]

Define
\[
 \widehat z=z+\sum_{i=1}^mX_iY_i.
\]
The grading relations give
\[
 [\widehat z,X_i]=[\widehat z,Y_i]=0,
 \qquad [\widehat z,c]=-2c.
\]
The elements $z$ and $X_iY_i$ commute.  In a Weyl pair with $[X,Y]=1$,
writing $D=XY$, one has
\[
 X^rY^r=D(D+1)\cdots(D+r-1),
\]
so the polynomial identity
$\prod_{a\in\Fp}(T+a)=T^p-T$ yields
\[
 (XY)^p-XY=X^pY^p.
\]
Consequently $\widehat z^p-\widehat z\in k1$.  Set
$h_0=-\frac12\widehat z$ and choose $\lambda\in k$ such that
$h=h_0+\lambda1$ satisfies $h^p-h=0$.  With $u=a^{-1}c$,
\[
 [h,u]=u,\qquad h^p-h=0,\qquad u^p=1.
\]
By \cref{lem:toral-shift}, $h,u$ generate a unital $\Mat_p(k)$ commuting
with $M_-$.  Multiplication gives a unital map
\[
 \Mat_{p^m}(k)\otimes\Mat_p(k)\longrightarrow\Achi.
\]
Its domain is simple, so the map is injective.
\end{proof}

\subsection{Vector negative support and symplectic automorphisms}

\begin{lemma}[Restricted symplectic substitutions]\label{lem:sp-auto}
For $\sigma\in\operatorname{Sp}(V)$ define
\[
 (T_\sigma f)(v,x_n)=f(\sigma^{-1}v,x_n).
\]
Then
\[
 D_K(f)\longmapsto D_K(T_\sigma f)
\]
is a restricted automorphism $\phi_\sigma$ of $\kk$, and it induces
\[
 u_\chi(\kk)\cong u_{\chi\circ\phi_\sigma^{-1}}(\kk).
\]
\end{lemma}

\begin{proof}
Linear substitution preserves the height-one truncation because
\[
 \Bigl(\sum_i a_ix_i\Bigr)^p=\sum_i a_i^px_i^p=0.
\]
The map $T_\sigma$ fixes $x_n$, commutes with $E$ and $\partial_{x_n}$,
and symplecticity gives
\[
 \{T_\sigma f,T_\sigma g\}_{\mathrm{sp}}
 =T_\sigma\{f,g\}_{\mathrm{sp}}.
\]
Equation \eqref{eq:contact-bracket} therefore gives
\[
 \{T_\sigma f,T_\sigma g\}_K
 =T_\sigma\{f,g\}_K.
\]
In the derivation realization this is conjugation by $T_\sigma$, which
commutes with $p$th powers.  The asserted isomorphism of reduced enveloping
algebras follows.
\end{proof}

\begin{proposition}\label{prop:vector-negative}
Assume $\chi(c)=0$ and $\chi|_V\ne0$.  Then $\Achi$ contains a unital
copy of $\Mat_{p^m}(k)$.
\end{proposition}

\begin{proof}
Using \cref{lem:sp-auto} and the transitivity of $\operatorname{Sp}(V)$ on
nonzero vectors, normalize
\[
 \chi(u_1)=1,\quad \chi(v_1)=0,\quad
 \chi(u_j)=\chi(v_j)=0\quad(j\ge2).
\]
Then $u_1^p=1$.  For $2\le j\le m$ choose
$A_j\in\mathfrak{sp}(V)$ with
\[
 A_j(v_j)=u_1,\qquad A_j(v_1)=u_j,
\]
and zero on the other basis vectors.  Then
\[
 [A_i,A_j]=0,\qquad A_j^{[p]}=0,
 \qquad [A_i,v_j]=\delta_{ij}u_1,
 \qquad [A_j,u_1]=0.
\]
Put
\[
 X_j=u_1^{-1}A_j,
 \qquad P_j=X_j-\chi(A_j)1,
 \qquad Q_j=v_j.
\]
Then
\[
 P_j^p=Q_j^p=0,
 \qquad [P_i,Q_j]=\delta_{ij},
\]
and the different pairs commute.  They generate
$\Mat_{p^{m-1}}(k)$.

Let $h_1\in\mathfrak{sp}(V)$ act by $+1$ on $u_1$, by $-1$ on $v_1$,
and trivially on the remaining symplectic basis.  Then $h_1^{[p]}=h_1$,
$[h_1,u_1]=u_1$, and $h_1$ commutes with the preceding $m-1$ Weyl pairs.
Choose $\lambda$ so that $h=h_1+\lambda1$ satisfies $h^p-h=0$.  By
\cref{lem:toral-shift}, $h,u_1$ generate a further unital $\Mat_p(k)$
commuting with the preceding factor.  Their tensor product embeds because it
is simple.
\end{proof}

\section{Positive-leading characters}\label{sec:positive}

Let $d>0$ be the highest degree on which $\chi$ is nonzero and write
$\lambda=\chi|_{K_d}$.  A homogeneous element of $\kk$ of degree $r>0$
has restricted $p$-power in degree $pr$.  By \cref{lem:leading-unit},
every homogeneous $u\in K_d$ with $\lambda(u)\ne0$ is invertible in
$\Achi$.  For any $Y$ with $Y^{[p]}=0$ set
\[
 \bar Y=Y-\chi(Y)1.
\]
Then $\ad\bar Y=\ad Y$ and the reduced relation gives $\bar Y^p=0$.

\begin{proposition}[Nonresonant toral shift]\label{prop:positive-nonres}
If $p\nmid d$, then $\Achi$ contains a unital $\Mat_p(k)$.
\end{proposition}

\begin{proof}
Choose $u\in K_d$ with $\chi(u)\ne0$.
Put $h=d^{-1}z$.  Since $z^{[p]}=z$ and $d\in\Fp^\times$,
\[
 h^{[p]}=d^{-p}z^{[p]}=d^{-1}z=h,
\]
and $[h,u]=u$.  Apply \cref{lem:leading-weight}.
\end{proof}

\begin{lemma}[Saturated-support square-zero lemma]\label{lem:saturated}
Let
\[
 \varnothing\ne S\subseteq\{1,\ldots,m\}
\]
and
\[
 \psi_S=\prod_{i\in S}u_i^{(p-1)}v_i^{(p-1)}.
\]
A coefficient is called $S$-saturated-first-order if in every plane of $S$
all exponents are $p-1$ except possibly one exponent $p-2$.  The product of
two such coefficients is zero in the truncated divided-power algebra, and
the product of one such coefficient with a first derivative of another is
zero.  Consequently any derivation all of whose coefficients are
$S$-saturated-first-order satisfies $X^2=0$.
\end{lemma}

\begin{proof}
The smallest possible exponent sum on a common coordinate is
$(p-2)+(p-1)=2p-3\ge p$.  If a first derivative creates exponent $p-3$,
the worst same-coordinate sum is $2p-5\ge p$, while if that coordinate does
not truncate, another coordinate in the same saturated plane has sum at
least $2p-3$.  Hence all products occurring in $X^2$ vanish.
\end{proof}

\begin{proposition}[Even resonant zero-weight case]\label{prop:even-resonant}
Suppose $d=\ell p$ with $\ell=2s$, and suppose the leading functional is supported
on the joint zero-weight subspace of the diagonal symplectic torus.  Then
$\Achi$ contains a unital $\Mat_p(k)$.
\end{proposition}

\begin{proof}
For a monomial $u^Av^Bx_n^{(c)}$, the $h_i$-weight is
$B_i-A_i$.  Joint zero weight therefore forces $A_i=B_i=a_i$ for every
$i$.  Its Contact weight is then
\[
 2\sum_i a_i+2c.
\]
Since the functions representing $K_{\ell p}$ have Contact weight
$\ell p+2=2sp+2$, one obtains
\[
 2\sum_i a_i+2c=2sp+2,
 \qquad\text{hence}\qquad
 \sum_i a_i+c=sp+1.
\]
Thus the zero-weight leading functions have monomial basis
\[
 f_{A,c}=\prod_{i=1}^m u_i^{(a_i)}v_i^{(a_i)}x_n^{(c)},
 \qquad
 \sum_i a_i+c=sp+1,
 \qquad 0\le a_i,c\le p-1.
\]
A pure $x_n$-power is impossible because $\ell p+2>2(p-1)$.  Hence any
monomial on which $\lambda$ is nonzero uses at least one symplectic plane.

If some $1\le a_j\le p-2$, apply the symplectic transvection
$u_j'=u_j+v_j$, $v_j'=v_j$.  Replacing the plane factor by
$(u_j')^{(a_j+1)}(v_j')^{(a_j-1)}$ produces a vector of nonzero weight for
$D_K(u_j'v_j')$.  Indeed,
\[
 (u_j+v_j)^{(a_j+1)}v_j^{(a_j-1)}
 =\sum_{i=0}^{a_j+1}
 \binom{2a_j-i}{a_j-1}u_j^{(i)}v_j^{(2a_j-i)}.
\]
The old zero-weight monomial occurs only for $i=a_j$, and its coefficient is
\[
 \binom{a_j}{a_j-1}=a_j\in k^\times.
\]
Hence the leading functional remains nonzero on the transformed
nonzero-weight vector, and the toral-shift argument applies.

Thus $a_i\in\{0,p-1\}$ may be assumed for every $i$.  Let
$S=\{i:a_i=p-1\}$, which is nonempty, and put
$\psi_S=\prod_{i\in S}u_i^{(p-1)}v_i^{(p-1)}$.

If $c\le p-2$, set
\begin{align*}
 f&=x_n^{(c)}\psi_S,
 &F&=x_n^{(c+1)}\psi_S,
 &Y&=D_K(1),\\
 \bar Y&=Y-\chi(Y)1,
 &R&=D_K(F).
\end{align*}
Then $\bar Y^p=0$ and $\delta:=\ad\bar Y=\ad Y$ satisfies
$\delta^p=0$, since $(\ad\bar Y)^p=\ad(\bar Y^p)$ in characteristic
$p$.  Since $R$ lies strictly above the leading degree $d$, one has
$\chi(R)=0$.  In the derivation realization the restricted $p$-map is the
$p$th composition; \cref{lem:saturated} gives $R^2=0$, hence
$R^{[p]}=0$.  The reduced relation therefore gives $R^p=0$ in $\Achi$.  Moreover
\[
 \delta(R)=[Y,D_K(F)]=2D_K(f),
\]
and the right-hand side is a unit because $\lambda(D_K(f))\ne0$.

To verify the Hasse hypothesis, every
$\delta^j(R)$ is a scalar multiple of a derivation
$D_K(x_n^{(c+1-j)}\psi_S)$ (or is zero).  Along this Hasse string only the
$x_n$-exponent changes; the symplectic exponents in $\psi_S$ remain saturated.
Every residue has the coefficient pattern of \cref{lem:saturated}.  Its truncation estimate therefore applies to any two
Hasse residues: every term in their Contact bracket contains either a product
of two $S$-saturated coefficients or a product of one such coefficient with a
first derivative of another, and hence vanishes.  Hence the finite algebra
\[
 H=k[R,\delta R,\ldots,\delta^{p-1}R]
\]
is commutative.  In the Hasse expansion of
$\delta^{p-1}(R^{p-1})$, a word is indexed by
\[
 a_1+\cdots+a_{p-1}=p-1,
 \qquad a_i\ge0,
\]
and the principal word is the unique one with
$a_1=\cdots=a_{p-1}=1$.  Every nonprincipal word therefore has some
$a_i=0$ and hence contains an undifferentiated factor $R$.  Thus every
nonprincipal word lies in the ideal $(R)\subset H$.  Since $H$ is
commutative and $R^p=0$,
\[
 (R)^p=(R^p)=0.
\]
The element $\delta(R)$ lies in $H$ and is invertible in $\Achi$.
Multiplication by $\delta(R)$ is therefore injective, hence bijective, on the
finite-dimensional algebra $H$; in particular, $\delta(R)^{-1}\in H$.
Normalization by the principal term preserves $(R)$, so the normalized error
is nilpotent.  The assertion follows from \cref{lem:hasse-slice}.

Define the full-top Lie degree
\[
 \Delta_K:=2(m+1)(p-1)-2.
\]
At $c=p-2$ with all symplectic planes saturated,
\[
 \psi_{\{1,\ldots,m\}}
 =\prod_{i=1}^m u_i^{(p-1)}v_i^{(p-1)},
 \qquad
 x_n^{(p-1)}\psi_{\{1,\ldots,m\}}=x^{(\tau)},
\]
where $\tau=(p-1,\ldots,p-1)$.  Thus the required antiderivative in $x_n$ is
the full-top function.  The corresponding one-hole monomial has degree
$\Delta_K-2$, and the resonant congruence is
\[
 \Delta_K-2\equiv-2(m+3)\equiv0\pmod p,
\]
whereas
\[
 \Delta_K\equiv-2(m+2)\not\equiv0\pmod p.
\]
Since $\Delta_K\not\equiv0\pmod p$, the full-top function is not deleted on
passage to $K(2m+1;\mathbf1)^{(1)}$, so the required antiderivative lies in
the derived Contact algebra.

Suppose now that $c=p-1$.  If $S\ne\{1,\ldots,m\}$, choose
$r\notin S$ and put
\[
 \zeta_r=x_n-u_rv_r,
 \qquad f_0=\zeta_r^{(p-1)}\psi_S.
\]
We may assume that the leading functional vanishes on every zero-weight
monomial treated above.  In the expansion of $f_0$, all terms other than the
original saturated
$c=p-1$ monomial belong to those earlier subcases.  Therefore their
$\lambda$-values vanish, whereas the original term has nonzero
$\lambda$-value, and hence
\[
 \lambda(D_K(f_0))\ne0.
\]
Set
\[
 F=v_rf_0,
 \quad Y=D_K(u_r),
 \quad \bar Y=Y-\chi(Y)1,
 \quad X=D_K(F),
 \quad Z=D_K(f_0).
\]
Because $f_0$ depends on $u_r,v_r$ only through $\zeta_r$, one has
\[
 D_{u_r}f_0=-v_rD_{x_n}f_0,
 \qquad
 D_{v_r}f_0=-u_rD_{x_n}f_0.
\]
Substitution into \eqref{eq:contact-bracket} gives
\[
 \{u_r,v_rf_0\}_K=f_0,
 \qquad
 \{u_r,f_0\}_K=0,
 \qquad
 \{v_rf_0,f_0\}_K=0.
\]
Consequently
\[
 [\bar Y,X]=Z,
 \qquad [\bar Y,Z]=0,
 \qquad [X,Z]=0.
\]
The scalar centering gives $\bar Y^p=0$.  The element $X$ has degree strictly larger than $d$, so $\chi(X)=0$.  As
above, the derivation realization identifies the restricted $p$-map with the
$p$th composition, and \cref{lem:saturated} gives $X^2=0$, hence
$X^{[p]}=0$.  The reduced relation therefore yields
\[
 X^p=X^{[p]}+\chi(X)^p1=0.
\]
Since $Z$ is a leading unit and commutes with $\bar Y$ and $X$, the elements
\[
 P=Z^{-1}\bar Y,
 \qquad Q=X
\]
satisfy
\[
 [P,Q]=Z^{-1}[\bar Y,X]=1,
\]
\[
 P^p=Z^{-p}\bar Y^p=0,
 \qquad Q^p=X^p=0.
\]
Thus \cref{lem:weyl-matrix} applies.

The remaining monomial, with $S=\{1,\ldots,m\}$ and $c=p-1$, is
\[
 x_n^{(p-1)}\prod_{i=1}^m u_i^{(p-1)}v_i^{(p-1)}=x^{(\tau)},
\]
so it is the full-top function itself and has Lie degree $\Delta_K$.
Resonance would force
\[
 \Delta_K\equiv-2(m+2)\equiv0\pmod p.
\]
When $2m+4\equiv0\pmod p$, however, this is precisely the basis element
deleted on passage to $K(2m+1;\mathbf1)^{(1)}$; see
\cite[Section~2]{Mygind2014}.
\end{proof}

\begin{corollary}\label{cor:positive-leading}
If the highest nonzero homogeneous component of $\chi$ has positive degree,
then $u_\chi(\kk)$ contains a unital copy of $\Mat_p(k)$.
\end{corollary}

\begin{proof}
If $p\nmid d$, apply \cref{prop:positive-nonres}.
Suppose $d=\ell p$.  On a monomial $u^Av^Bx_n^{(c)}$, the weight for
$h_i=D_K(u_iv_i)$ is $B_i-A_i$.  A joint zero-weight monomial has
$A_i=B_i$, hence Contact weight $2\sum_iA_i+2c$.  The functions in
$K_{\ell p}$ have Contact weight $\ell p+2$, which is odd when $\ell$ is odd, so
the joint zero-weight space vanishes for odd $\ell$.  Thus if $\ell$ is odd,
or if $\ell$ is even and $\lambda$ does not vanish on the nonzero
torus-weight space, $\lambda$ detects a nonzero torus-weight vector $u$.
For some toral $h$, $[h,u]=cu$ with $c\in\Fp^\times$.  Apply
\cref{lem:leading-weight}.
If $\ell$ is even and $\lambda$ is supported on the joint zero-weight
subspace, apply \cref{prop:even-resonant}.
\end{proof}

\section{Characters supported on \texorpdfstring{$K_0$}{K0}}\label{sec:k0}

Assume now
\[
 \chi(K_j)=0\quad(j\ne0),
 \qquad 0\ne\chi|_{K_0}.
\]
Recall
\[
 K_0\cong\mathfrak{sp}(V)\oplus kz.
\]

\begin{proposition}[Nonzero symplectic degree-zero part]\label{prop:k0-sp}
If $\chi|_{\mathfrak{sp}(V)}\ne0$, then $\Achi$ contains a unital
$\Mat_p(k)$.
\end{proposition}

\begin{proof}
Under $\mathfrak{sp}(V)\cong S^2V$, pure squares span $S^2V$ because
$2\in k^\times$.  Thus for some $0\ne w\in V$,
\[
 \chi(D_K(w^{(2)}))\ne0.
\]
By \cref{lem:sp-auto}, normalize $w=u_1$.  Put
\[
 e=D_K(u_1^{(2)}),\qquad h=D_K(u_1v_1),
\]
so that
\[
 [h,e]=de,
 \qquad e^{[p]}=0,
 \qquad h^{[p]}=h,
 \qquad d\in\Fp^\times.
\]
Put $b=\chi(e)\ne0$ and $u=b^{-1}e$, so $u^p=1$.  Choose
$\lambda$ with $\lambda^p-\lambda=\chi(h)^p$ and put
$H=d^{-1}(h-\lambda1)$.  Then
\[
 [H,u]=u,
 \qquad H^p-H=0,
 \qquad u^p=1.
\]
Apply \cref{lem:toral-shift}.
\end{proof}

\begin{proposition}[Pure grading-line degree-zero part]\label{prop:k0-z}
If $\chi|_{\mathfrak{sp}(V)}=0$ and $a=\chi(z)\ne0$, then $\Achi$
contains a unital $\Mat_p(k)$.
\end{proposition}

\begin{proof}
Let
\[
 P=D_K(u_1),
 \qquad h_1=D_K(u_1v_1),
 \qquad H=z-h_1.
\]
Then
\[
 [P,H]=0,
 \qquad H^{[p]}=H,
 \qquad \chi(H)=a.
\]
Define
\[
 X=D_K(v_1x_n-2u_1v_1^{(2)}).
\]
The Contact bracket gives
\[
 [P,X]=H.
\]
In $\Achi$,
\[
 H^p-H=a^p1,
\]
so
\[
 H^{-1}=a^{-p}(H^{p-1}-1).
\]
Set $Q=XH^{-1}$.  Since $[P,H]=0$,
\[
 [P,Q]=1,
 \qquad P^p=0.
\]
Apply \cref{lem:one-nilpotent} to $P,Q$.
\end{proof}

\begin{proof}[Proof of \cref{thm:first-matrix}]
The central and vector negative-support cases are
Propositions~\ref{prop:central-negative} and~\ref{prop:vector-negative}.
The positive-leading case is Corollary~\ref{cor:positive-leading}.  In the
remaining case $\chi$ is supported on $K_0$; Propositions~\ref{prop:k0-sp}
and~\ref{prop:k0-z} apply.
\end{proof}

\section{The Witt--Jacobson, Special and Hamiltonian families}\label{sec:wsh}

The following constructions give the first matrix factor for
$W(n)$, $S(n)^{(1)}$ and $H(2m)^{(2)}$.

\subsection{The Witt--Jacobson algebra}

Let $A(n)=k[x_1,\ldots,x_n]/(x_1^p,\ldots,x_n^p)$ and
$\mathfrak w=W(n)=\operatorname{Der}A(n)$, with the divided-monomial basis
$x^\alpha\partial_i$, $0\le\alpha_j<p$, and grading
$\mathfrak w_r=\spn\{x^\alpha\partial_i:|\alpha|=r+1\}$,
$-1\le r\le N=n(p-1)-1$; see \cite[Section~2]{Mygind2014}.  A homogeneous
derivation of degree $r$ has restricted $p$-power homogeneous of degree $pr$
or zero.  Put $h_j=x_j\partial_j$; then $h_j^{[p]}=h_j$ and
$[h_j,x^\alpha\partial_i]=(\alpha_j-\delta_{ij})x^\alpha\partial_i$.

\begin{theorem}[Witt--Jacobson first matrix factor]\label{thm:witt-first-matrix}
For every $n\ge1$ and every $0\ne\chi\in\mathfrak w^*$,
\[
 u_\chi(\mathfrak w)\supseteq\Mat_p(k)
\]
with the same unit.
\end{theorem}

\begin{proof}
\emph{Case $\chi(\mathfrak w_{-1})\ne0$.}  Choose $i$ with
$a=\chi(\partial_i)\ne0$.  Since $\partial_i^{[p]}=0$, in $u_\chi(\mathfrak
w)$ one has $\partial_i^p=a^p1$, so $U=a^{-1}\partial_i$ is a unit.  Put
$h=-x_i\partial_i$; then $h^{[p]}=h$ and $[h,U]=U$.  Apply
\cref{cor:toral-shift-strong}.

\emph{Case $\chi(\mathfrak w_{-1})=0$, $\chi(\mathfrak w_{>0})\ne0$.}  Let
$r>0$ be maximal with $\chi(\mathfrak w_r)\ne0$ and choose
$\Theta=x^\alpha\partial_i\in\mathfrak w_r$ with $\chi(\Theta)\ne0$.  The
weight $\lambda=\alpha-\varepsilon_i$ is nonzero, since $\lambda=0$ forces
$\alpha=\varepsilon_i$ and $r=0$.  Choose $j$ with
$c=\alpha_j-\delta_{ij}\ne0$.  Then $h_j^{[p]}=h_j$ and
$[h_j,\Theta]=c\Theta$.  Apply \cref{lem:leading-weight}.

\emph{Case $\chi(\mathfrak w_{-1})=\chi(\mathfrak w_{>0})=0$,
$\chi|_{\mathfrak w_0}\ne0$.}  If $n\ge2$, choose $a,b$ with
$\chi(x_a\partial_b)\ne0$ and $s\ne a$.  The endomorphism
$\phi_c(x_b)=x_b+cx_sx_a$, $\phi_c(x_j)=x_j$ $(j\ne b)$, is an automorphism
of $A(n)$ fixing $\mathfrak m/\mathfrak m^2$, and conjugation
$g_c(D)=\phi_cD\phi_c^{-1}$ is a restricted automorphism of $\mathfrak w$.
A direct computation gives
$g_c^{-1}(\partial_s)=\partial_s+cx_a\partial_b+Y$ with $Y\in\mathfrak
w_{\ge1}$, so
\[
 (g_c\cdot\chi)(\partial_s)=\chi(g_c^{-1}\partial_s)=c\,\chi(x_a\partial_b)\ne0
\]
for $c\ne0$, and $g_c\cdot\chi$ has nonzero negative support.  If $n=1$,
take $\phi_c(x)=x+cx^2$; then
$\phi_c^{-1}\partial\phi_c=(1+2cx+\cdots)\partial$, and
$(g_c\cdot\chi)(\partial)=2c\,\chi(x\partial)\ne0$ since $p\ge5$.  In either
case the restricted automorphism identifies $u_\chi(\mathfrak w)$ with
$u_{g_c\cdot\chi}(\mathfrak w)$, which falls under the first case.
\end{proof}

\begin{example}[A rank-two matrix factor]\label{ex:w2-negative}
Let $\mathfrak g=W(2)$ and take $\chi(\partial_1)=1$, with
$\chi$ zero on every other monomial basis vector. In $A=u_\chi(\mathfrak g)$ put
\[
 U=\partial_1,\qquad H=-x_1\partial_1,\qquad
 P=\partial_2,\qquad Q=(x_2\partial_1)U^{-1}.
\]
The reduced relations and the derivation brackets give
\[
 U^p=1,\quad H^p=H,\quad [H,U]=U,\qquad
 P^p=Q^p=0,\quad [P,Q]=1.
\]
The element $U$ commutes with $x_2\partial_1$, whose restricted power
and character value vanish. Also $[H,x_2\partial_1]=x_2\partial_1$
and $[H,U^{-1}]=-U^{-1}$, so the two pairs commute.
The two pairs give a unital $\Mat_{p^2}(k)$ and a corner $C$ with
\[
 A\cong\Mat_{p^2}(C),\qquad \dim C=p^{2p^2-4}.
\]
For $p=5$ this is a $25$ by $25$ matrix factor with corner dimension
$5^{46}$; for $p=7$ it is a $49$ by $49$ matrix factor with corner
dimension $7^{94}$. Every simple $A$-module has dimension divisible by $p^2$.
The form $\chi([\cdot,\cdot])$ pairs $\partial_1$ with $x_1\partial_1$
and $\partial_2$ with $x_2\partial_1$, and vanishes on all other pairs
of monomial basis vectors. Hence $\dim\mathfrak g^\chi=\dim\mathfrak g-4$,
and $p^2$ is the Kac--Weisfeiler bound at $\chi$.
\end{example}

\subsection{The Special algebra}

For $n\ge3$ let $D_{ij}(f)=\partial_j(f)\partial_i-\partial_i(f)\partial_j$
$(i\ne j)$ generate $\mathfrak s=S(n)^{(1)}\subset W(n)$, with
$\mathfrak s_0\cong\mathfrak{sl}_n$ and torus generated by
$x_a\partial_a-x_b\partial_b$, as in Mygind
\cite[Section~2]{Mygind2014}. Ma determined the graded bracket
$[\mathfrak s_2,\mathfrak s_{d-1}]$ over a general field
\cite{Ma2026Divergence}. The identities below are the root-graded forms
in positive characteristic. Write
$q=x_1\cdots x_n$ and $E_{ab}=x_a\partial_b$ $(a\ne b)$ for the root vectors
of $\mathfrak s_0$.  Expanding the definitions of $E_{ab}$ and $D_{ij}$
yields the equivariance identity
\begin{equation}\label{eq:special-equiv}
 [E_{ab},D_{ij}(f)]=D_{ij}(E_{ab}f)+\delta_{ia}D_{jb}(f)-\delta_{ja}D_{ib}(f).
\end{equation}

\begin{lemma}[Zero-weight classification]\label{lem:special-zero-weight}
The positive-degree torus-weight-zero subspace of $\mathfrak s$ is
spanned by the vectors
\[
 U_{ij,c}=D_{ij}(q^cx_ix_j),\qquad 1\le c\le p-2,\ i\ne j.
\]
\end{lemma}

\begin{proof}
The torus weight of $D_{ij}(x^\beta)$ is $\beta-\varepsilon_i-\varepsilon_j$
modulo the all-ones vector, so it vanishes exactly when
$\beta=c\mathbf1+\varepsilon_i+\varepsilon_j$ for some $c\in\Fp$ (with
exponents reduced into $[0,p-1]$).  For $c=0$ this is $x_ix_j$, which lies
in $\mathfrak s_0$; for $c=p-1$ both exponents $\beta_i,\beta_j$ vanish and
$D_{ij}(x^\beta)=0$; for $1\le c\le p-2$, $x^\beta=q^cx_ix_j$ and
$D_{ij}(q^cx_ix_j)=(c+1)q^c(x_i\partial_i-x_j\partial_j)\ne0$ lies in the
positive part.
\end{proof}

\begin{proposition}[Zero-weight root-span]\label{prop:special-root-span}
Every $U_{ij,c}$ lies in
\[
 \sum_{\alpha\ne0}[\mathfrak s_0(\alpha),\mathfrak s(-\alpha)].
\]
\end{proposition}

\begin{proof}
If $n\ge4$, choose $a,b\notin\{i,j\}$ distinct and set
$\gamma=c\mathbf1+\varepsilon_i+\varepsilon_j-\varepsilon_a+\varepsilon_b$,
which is legal since $1\le c\le p-2$.  With $V=D_{ij}(x^\gamma)$, one has
$\operatorname{wt}(V)=\varepsilon_b-\varepsilon_a\ne0$ and, since
$a,b\notin\{i,j\}$, \eqref{eq:special-equiv} gives
$[E_{ab},V]=(c+1)U_{ij,c}$.

If $n=3$, relabel $(i,j,k)=(1,2,3)$.  For $1\le c\le p-3$, put
$\gamma=c\mathbf1+\varepsilon_i+\varepsilon_j-\varepsilon_k+\varepsilon_i$;
\eqref{eq:special-equiv} gives
$[E_{ki},D_{ij}(x^\gamma)]=(c+2)U_{ij,c}$.  For the endpoint $c=p-2$, set
$U_{12}=D_{12}(q^cx_1x_2)$ and $U_{23}=D_{23}(q^cx_2x_3)$.  A direct
computation from \eqref{eq:special-equiv} and the definitions gives
$U_{13}=D_{13}(q^cx_1x_3)=U_{12}+U_{23}$; the two contributions along
$\partial_2$ cancel. Then
\[
 [E_{13},D_{12}(q^cx_2x_3)]=(c+1)U_{12}+U_{23},\qquad
 [E_{23},D_{12}(q^cx_1x_3)]=cU_{12}-U_{23}.
\]
In the basis $(U_{12},U_{23})$ these rows are $(-1,1)$ and $(-2,-1)$, with
determinant $3\ne0$ since $p\ge5$; hence $U_{12},U_{23}$, and therefore
$U_{13}$, are root images.
\end{proof}

\begin{lemma}[Tangent transport]\label{lem:tangent-transport}
Let $t\mapsto g_t$ be a polynomial family of restricted automorphisms of
a restricted Lie algebra, with $g_0=\mathrm{id}$, whose derivative at
$t=0$ is $\operatorname{ad}E$.  If $\chi([E,V])\ne0$, then
$\chi(g_t(V))\ne0$ for some $t\in k$.
\end{lemma}

\begin{proof}
The value $t\mapsto\chi(g_t(V))$ is a polynomial, with derivative
$\chi([E,V])$ at zero.  Since $k$ is infinite, the polynomial does not
vanish at every $t$.
\end{proof}

\begin{theorem}[Special first matrix factor]\label{thm:special-first-matrix}
For every $n\ge3$ and every $0\ne\chi\in\mathfrak s^*$,
\[
 u_\chi(\mathfrak s)\supseteq\Mat_p(k)
\]
with the same unit.
\end{theorem}

\begin{proof}
\emph{Negative support.} If $\chi(\partial_s)\ne0$, then $\partial_s$ gives
a unit as in \cref{thm:witt-first-matrix}, and $h=x_t\partial_t-x_s\partial_s$
($t\ne s$) satisfies $h^{[p]}=h$, $[h,\partial_s]=\partial_s$; apply
\cref{cor:toral-shift-strong}.

\emph{Positive-leading, nonzero weight.} If the top nonvanishing degree
$r>0$ has $\chi(\Theta)\ne0$ for some torus-weight vector
$\Theta\in\mathfrak s_r(\lambda)$, $\lambda\ne0$, choose a toral
$h\in\mathfrak s_0$ with $[h,\Theta]=\Theta$.  As a graded restricted
subalgebra of $W(n)$, $\mathfrak s$ satisfies the hypotheses of
\cref{lem:leading-weight}.  Apply that lemma.

\emph{Positive-leading, zero weight.} If the top support is torus-weight
zero, \cref{prop:special-root-span} gives $E\in\mathfrak s_0(\alpha)$ and
$V\in\mathfrak s_r(-\alpha)$ with $\chi([E,V])\ne0$.  Let $g_t$ be
conjugation by the linear substitution $x_b\mapsto x_b+tx_a$ for the
root vector $E=E_{ab}$.  Its determinant is one, so it preserves
$\mathfrak s$ and its grading, and its tangent action at zero is
$\operatorname{ad}E$.  By \cref{lem:tangent-transport},
$\chi(g_t(V))\ne0$ for some $t$.  The transported character detects the
nonzero-weight vector $V$ in degree $r$, and the preceding case applies.

\emph{Pure degree zero.} Suppose $\chi(\mathfrak s_{-1})=
\chi(\mathfrak s_{>0})=0$ and $\chi|_{\mathfrak s_0}\ne0$.
First arrange that $\chi(E_{ab})\ne0$ for some root vector.
If $\chi$ initially vanishes on every root space, choose
$h_{ij}=x_i\partial_i-x_j\partial_j$ with $\chi(h_{ij})\ne0$.
Conjugation by $x_i\mapsto x_i+t x_j$ sends
\[
 E_{ij}\longmapsto E_{ij}-t h_{ij}-t^2E_{ji}.
\]
For $t\ne0$ the transported character is nonzero on $E_{ij}$.

Now choose $s\notin\{a,b\}$ and use the polynomial substitution
\[
 \theta_t(x_b)=x_b-tx_sx_a,\qquad
 \theta_t(x_j)=x_j\quad(j\ne b).
\]
Its inverse is $\theta_{-t}$ and its Jacobian determinant is one.
Conjugation $g_t(D)=\theta_tD\theta_t^{-1}$ is therefore a restricted
automorphism of $\mathfrak s$. Evaluation on the coordinate generators gives
\[
 g_t(\partial_s)=\partial_s+tE_{ab},\qquad
 (\chi\circ g_t)(\partial_s)=t\chi(E_{ab})\ne0.
\]
The negative-support case applies.
\end{proof}

\subsection{The Hamiltonian algebra}

Let $D_H$ denote the Hamiltonian map in symplectic coordinates. The
Poisson bracket
\[
 \{f,g\}=\sum_k\bigl(\partial_{x_k}f\,\partial_{y_k}g
 -\partial_{y_k}f\,\partial_{x_k}g\bigr)
\]
is the action of $D_H(f)$ on $g$, and the map is normalized by
$[D_H(f),D_H(g)]=D_H(\{f,g\})$.
Then $\mathfrak h=H(2m)^{(2)}$ has $\mathfrak h_0\cong\mathfrak{sp}_{2m}$;
see \cite[Section~2]{Mygind2014}.
Put $t_i=D_H(x_iy_i)$. The bracket gives
$t_i=y_i\partial_{y_i}-x_i\partial_{x_i}$, hence
$t_i(x^\alpha)=(\alpha_{y_i}-\alpha_{x_i})x^\alpha$.
The normalization then gives
\begin{equation}\label{eq:ham-weight}
 [t_i,D_H(x^\alpha)]=(\alpha_{y_i}-\alpha_{x_i})D_H(x^\alpha).
\end{equation}
The coefficient $\alpha_{y_i}-\alpha_{x_i}$ is the $t_i$-weight of
$D_H(x^\alpha)$.

\begin{lemma}[Zero-weight classification]\label{lem:hamiltonian-zero-weight}
The positive-degree torus-weight-zero subspace of the simple derived
Hamiltonian algebra is spanned by
$D_H(f_{\mathbf a})$, $f_{\mathbf a}=\prod_i(x_iy_i)^{a_i}$, with
$\sum_i a_i\ge2$ and $\mathbf a\ne(p-1,\ldots,p-1)$.
\end{lemma}

\begin{proof}
Zero $t_i$-weight for exponents in $[0,p-1]$ forces $\alpha_{x_i}=\alpha_{y_i}$
for every $i$, giving the stated monomials; the constant tuple has zero
Hamiltonian vector field and the all-top tuple is excluded from the simple
derived algebra.
\end{proof}

The long- and short-root vectors below are themselves images under $D_H$:
\[
 E_i=x_i\partial_{y_i}=D_H\bigl(\tfrac12x_i^2\bigr),
 \qquad
 E_{ij}=x_i\partial_{x_j}-y_j\partial_{y_i}=D_H(-x_iy_j)\quad(i\ne j).
\]
Indeed $\{\tfrac12x_i^2,g\}=x_i\partial_{y_i}(g)$ and
$\{-x_iy_j,g\}=x_i\partial_{x_j}(g)-y_j\partial_{y_i}(g)$ for every $g$, by
direct computation from the Poisson bracket.
Consequently, for every $g$,
\begin{equation}\label{eq:ham-equiv}
 [E_i,D_H(g)]=D_H(\{\tfrac12x_i^2,g\})=D_H(E_i(g)),
 \qquad
 [E_{ij},D_H(g)]=D_H(E_{ij}(g)),
\end{equation}
by the identity $[D_H(f),D_H(g)]=D_H(\{f,g\})$.

\begin{proposition}[Zero-weight root-span]\label{prop:hamiltonian-root-span}
Every homogeneous $D_H(f_{\mathbf a})$ of degree $r>0$ in
\cref{lem:hamiltonian-zero-weight} lies in
$\sum_{\alpha\ne0}[\mathfrak h_0(\alpha),\mathfrak h_r(-\alpha)]$.
\end{proposition}

\begin{proof}
If $1\le a_i\le p-2$ for some $i$, let $g$
be obtained from $f_{\mathbf a}$ by replacing the $i$-th pair exponents
$(a_i,a_i)$ by $(a_i-1,a_i+1)$; then $E_i(g)=(a_i+1)f_{\mathbf a}$ by direct
evaluation, and \eqref{eq:ham-equiv} gives
$[E_i,D_H(g)]=(a_i+1)D_H(f_{\mathbf a})$ with $a_i+1\in k^\times$.

If every $a_s\in\{0,p-1\}$, then $\mathbf a\ne0$ and $\mathbf a\ne
(p-1,\ldots,p-1)$ force at least one index $i$ with $a_i=p-1$ and at least
one index $j\ne i$ with $a_j=0$.  Let $F$ be the product of the remaining paired factors,
and set $g=x_i^{p-2}x_jy_i^{p-1}F$.  Direct evaluation gives
\[
 E_{ij}(g)=x_i\partial_{x_j}(g)-y_j\partial_{y_i}(g)
 =f_{\mathbf a}-(p-1)x_i^{p-2}y_i^{p-2}x_jy_jF
 =f_{\mathbf a}+h,
\]
where $h=x_i^{p-2}y_i^{p-2}x_jy_jF$, using $-(p-1)\equiv1\pmod p$; and
\eqref{eq:ham-equiv} gives $[E_{ij},D_H(g)]=D_H(f_{\mathbf a}+h)$.  The
monomial $h$ has $j$-pair exponent $1$, so with $g_2$ obtained from $h$ by
replacing the $j$-pair exponents $(1,1)$ by $(0,2)$, direct evaluation
gives $E_j(g_2)=x_j\partial_{y_j}(g_2)=2x_jy_jx_i^{p-2}y_i^{p-2}F=2h$, and
\eqref{eq:ham-equiv} gives $[E_j,D_H(g_2)]=2D_H(h)$.  Hence
\[
 D_H(f_{\mathbf a})
 =[E_{ij},D_H(g)]-\tfrac12[E_j,D_H(g_2)]
 =D_H(f_{\mathbf a}+h)-D_H(h).
\]
Every term on the right has a nonzero-weight source.
\end{proof}

\begin{theorem}[Hamiltonian first matrix factor]\label{thm:hamiltonian-first-matrix}
For every $m\ge1$ and every $0\ne\chi\in\mathfrak h^*$,
\[
 u_\chi(\mathfrak h)\supseteq\Mat_p(k)
\]
with the same unit.
\end{theorem}

\begin{proof}
\emph{Negative support.} If $\chi(\partial_s)\ne0$ for $\partial_s\in
\{\partial_{x_i},\partial_{y_i}\}$, take the appropriate sign of
$t_i=D_H(x_iy_i)$ so that $[\pm t_i,\partial_s]=\partial_s$, and apply
\cref{cor:toral-shift-strong} to the resulting unit as in
\cref{thm:witt-first-matrix}.

\emph{Positive-leading, nonzero weight.} Let $r>0$ be maximal with
$\chi(\mathfrak h_r)\ne0$, and let $\Theta\in\mathfrak h_r$ be a
torus-weight vector with $\chi(\Theta)\ne0$ and weight $\lambda\ne0$.
Some component $c=\lambda_i$ lies in $\Fp^\times$.  Put
$t_i=D_H(x_iy_i)$ and $h=c^{-1}t_i$.  By \eqref{eq:ham-weight},
$[t_i,\Theta]=c\Theta$.  Since $t_i$ acts diagonally by weights in
$\Fp$, one has $t_i^{[p]}=t_i$, and the restricted $p$-map gives
$h^{[p]}=c^{-p}t_i=c^{-1}t_i=h$.  Thus $[h,\Theta]=\Theta$.
As a graded restricted subalgebra of $W(2m)$, $\mathfrak h$ satisfies
the hypotheses of \cref{lem:leading-weight}.  Apply that lemma.

\emph{Positive-leading, zero weight.}
By \cref{prop:hamiltonian-root-span}, there are $E,V$ with
$\chi([E,V])\ne0$, where $E$ is $E_i$ or $E_{ij}$.
If $E=E_i$, conjugate by $y_i\mapsto y_i+tx_i$.
If $E=E_{ij}$, conjugate by
\[
 x_j\mapsto x_j+tx_i,\qquad y_i\mapsto y_i-ty_j.
\]
The remaining variables stay fixed.  Both substitutions preserve the
symplectic form and have inverse given by $-t$, and the tangent action
at $t=0$ is $\operatorname{ad}E$.  \Cref{lem:tangent-transport} gives
$\chi(g_t(V))\ne0$ for some $t$.  The transported character falls under
the nonzero-weight case.

\emph{Pure degree zero.}
The degree-zero component contains the fields
\begin{equation}\label{eq:ham-zero-basis}
\begin{aligned}
&t_i=D_H(x_iy_i),
&&D_H(x_a^2),\ D_H(y_a^2),\\
&D_H(x_ax_b),\ D_H(y_ay_b)\ (a<b),
&&D_H(x_ay_b)\ (a\ne b).
\end{aligned}
\end{equation}
The coefficient of $\partial_{x_i}$ in $t_i$ is $-x_i$, so the $t_i$ are
linearly independent.  By \eqref{eq:ham-weight} the remaining vectors have
distinct torus weights for $p\ge5$, and none of those weights is zero, so
they are linearly independent from each other and from the $t_i$.  The list
has cardinality $m(2m+1)=\dim\mathfrak h_0$, and is therefore a basis.

If $\chi$ vanishes on every root vector, choose $i$ with $\chi(t_i)\ne0$.
Put $e_i=x_i\partial_{y_i}$ and $f_i=-y_i\partial_{x_i}=D_H(\tfrac12 y_i^2)$.
The substitution $\rho_s:y_i\mapsto y_i+sx_i$, fixing every other variable,
preserves the symplectic form and has inverse $\rho_{-s}$.  Conjugation gives
$g_s(f_i)=f_i+st_i+s^2e_i$.
Since $\chi$ vanishes on $f_i$ and on $e_i$, and $D_H(y_i^2)=2f_i$,
\[
 (\chi\circ g_s)(D_H(y_i^2))=2s\chi(t_i)\ne0\qquad(s\ne0).
\]
Replace $\chi$ by $\chi\circ g_s$.
After at most this replacement, $\chi(Z)\ne0$ for some root vector $Z$
in \eqref{eq:ham-zero-basis}.
Let $\varphi_a$ be the substitution $x_a\mapsto -y_a$, $y_a\mapsto x_a$,
fixing the other variables.  It preserves the symplectic form, and its
inverse is $x_a\mapsto y_a$, $y_a\mapsto -x_a$.  Write $g_a$ for conjugation
by $\varphi_a$, and $g_{ab}$ for conjugation by $\varphi_a\circ\varphi_b$.
For $b\ne a$,
\begin{align*}
 g_a\bigl(D_H(x_a^2)\bigr)&=D_H(y_a^2),\\
 g_a\bigl(D_H(x_ax_b)\bigr)&=-D_H(y_ax_b),\\
 g_b\bigl(D_H(x_ax_b)\bigr)&=-D_H(x_ay_b),\\
 g_{ab}\bigl(D_H(x_ax_b)\bigr)&=D_H(y_ay_b).
\end{align*}
The identity fixes $D_H(x_a^2)$ and $D_H(x_ax_b)$.  As $a$ and $b$ vary,
the four images on the right are scalar multiples of the remaining root
vectors in \eqref{eq:ham-zero-basis}.  Hence some substitution $g$ among
the identity, the $g_a$, and the $g_{ab}$ satisfies $g(R)=\lambda Z$ with
$\lambda\in k^\times$, where $R=D_H(x_i^2)$ or $R=D_H(x_ix_j)$ for some
$i\ne j$.  Replace $\chi$ by $\chi\circ g$.  Then $\chi(R)\ne0$.
Set $X=D_H(y_i)=-\partial_{x_i}$.
In the long-root case use
\[
 \sigma_t:y_i\mapsto y_i+tx_i^2;
\]
in the short-root case use
\[
 \sigma_t:y_i\mapsto y_i+tx_ix_j,\qquad
 y_j\mapsto y_j+\tfrac t2x_i^2.
\]
All other variables are fixed. These polynomial substitutions preserve
the symplectic form and have inverse $\sigma_{-t}$.
Conjugation gives, in both cases,
\[
 g_t(X)=X+tR,\qquad (\chi\circ g_t)(X)=t\chi(R)\ne0.
\]
The transported character has nonzero negative support.
\end{proof}

\begin{example}[A rank-two Hamiltonian matrix factor]\label{ex:h2-negative}
Let $\mathfrak h=H(2)^{(2)}$ in coordinates $x,y$, so that
$D_H(x)=\partial_y$, $D_H(y)=-\partial_x$, $D_H(xy)=y\partial_y-x\partial_x$
and $D_H(x^2)=2x\partial_y$. Take $\chi(D_H(x))=1$, with $\chi$ zero on
every other monomial basis vector. In $A=u_\chi(\mathfrak h)$ put
\[
 U=\partial_y,\qquad P=\partial_x,\qquad Q=(x\partial_y)U^{-1},\qquad
 H_0=x\partial_x-y\partial_y .
\]
Then $U^p=1$ and $P^p=0$. The elements $x\partial_y$ and $U$ commute and
$(x\partial_y)^p=0$, so $Q^p=0$; moreover $[P,Q]=\partial_yU^{-1}=1$
and $U$ commutes with $P$ and $Q$. By \eqref{eq:ham-weight},
\[
 [H_0,U]=U,\qquad [H_0,P]=-P,\qquad [H_0,Q]=Q .
\]
Unlike Example~\ref{ex:w2-negative}, $H_0$ does not commute with $P$ and $Q$, and the
torus of $\mathfrak h_0\cong\mathfrak{sl}_2$ is one-dimensional. Put
\[
 H=H_0-QP .
\]
Since $[P,Q]=1$ and $P^p=Q^p=0$, the element $QP$ is diagonalizable with
eigenvalues in $\Fp$, so $(QP)^p=QP$. It commutes with $H_0$, hence
$H^p=H$. Also $[QP,U]=0$, $[QP,P]=-P$ and $[QP,Q]=Q$, so
\[
 [H,U]=U,\qquad [H,P]=[H,Q]=0 .
\]
The pairs $(U,H)$ and $(P,Q)$ commute and each generates a unital copy of
$\Mat_p(k)$. Hence
\[
 A\cong\Mat_{p^2}(C),\qquad \dim C=p^{p^2-6}.
\]
For $p=5$ the corner has dimension $5^{19}$; for $p=7$ it has dimension
$7^{43}$. The form $\chi([\cdot,\cdot])$ pairs $D_H(x)$ with $D_H(xy)$
and $D_H(x^2)$ with $D_H(y)$, and vanishes on all other pairs of monomial
basis vectors, since $\{x^ay^b,x^cy^d\}$ has a nonzero coefficient of $x$
only for $a+c=2$ and $b+d=1$. Hence $\dim\mathfrak h^\chi=\dim\mathfrak h-4$,
and the factor $p^2$ dividing every simple $A$-module dimension is the
Kac--Weisfeiler bound at $\chi$.
\end{example}

\section{Peirce systems, symmetry and Morita equivalence}\label{sec:peirce}

For a finite-dimensional algebra $C$, write $\Soc(C)$ for the socle of the
left regular $C$-module.  For symmetric algebras this is also the right
regular socle.

\begin{lemma}[Peirce reconstruction]\label{lem:peirce}
Let $M=\spn\{E_{ij}\}_{0\le i,j<n}\cong\Mat_n(k)$ be a unital matrix
subalgebra of a unital algebra $A$, and let $e=E_{00}$.  Then
\[
 A\cong\Mat_n(eAe)
\]
via the coefficient maps
\[
 \gamma_{ij}(a)=E_{0i}aE_{j0}.
\]
In particular $e$ is full and $A$ is Morita equivalent to $eAe$.
\end{lemma}

\begin{proof}
The matrix-unit identities give
\[
 a=\sum_{i,j}E_{i0}\gamma_{ij}(a)E_{0j},
\]
and multiplication of coefficients is matrix multiplication.  Fullness
follows from $1=\sum_iE_{i0}eE_{0i}$.
\end{proof}

\begin{proposition}[Symmetry]\label{prop:symmetric}\label{prop:symmetric-general}
For $\mathfrak g\in\{W(n),S(n)^{(1)},H(2m)^{(2)},\kk\}$ and every
$\chi$, the algebra $u_\chi(\mathfrak g)$ and every full Peirce corner
reached from it are symmetric.
\end{proposition}

\begin{proof}
The reduced enveloping algebra of a finite-dimensional restricted Lie algebra
is Frobenius, with Nakayama automorphism determined on the Lie algebra by the
modular trace character $x\mapsto\operatorname{tr}(\ad x)$; see
\cite[Chapter~V, Section~4]{StradeFarnsteiner1988}.  Since $\mathfrak g$ is simple and nonabelian,
$[\mathfrak g,\mathfrak g]=\mathfrak g$. The modular trace character
therefore vanishes, so $u_\chi(\mathfrak g)$ is symmetric.  Symmetry is Morita invariant, so every full matrix
corner is symmetric.
\end{proof}

\begin{proposition}[Morita equivalence]\label{prop:morita-transport}
Suppose
\[
 A_j\cong\Mat_p(A_{j+1})\qquad(0\le j<r).
\]
Then $A_0\cong\Mat_{p^r}(A_r)$, and the successive corner functors yield
exact Morita equivalences.  For
corresponding modules,
\[
 \Ext_{A_0}^q(L,L')\cong\Ext_{A_r}^q(L_r,L_r')
 \qquad(q\ge0).
\]
They preserve blocks and projective covers, and
$\dim L=p^r\dim L_r$.
\end{proposition}

\begin{proof}
Composing the $r$ isomorphisms gives
\[
 A_0\cong\Mat_p(\Mat_p(\cdots\Mat_p(A_r)))=\Mat_{p^r}(A_r).
\]
For a full
matrix algebra $\Mat_N(B)$ the corner functor $M\mapsto eM$ attached to a
diagonal matrix unit $e=E_{11}$ is an equivalence
\[
 \Mat_N(B)\text{-mod}\longrightarrow B\text{-mod}
\]
with quasi-inverse $V\mapsto\Mat_N(B)e\otimes_BV$, since
$\Mat_N(B)e$ is a progenerator and $e\Mat_N(B)e=B$.  Both functors are exact, so the equivalence carries
$\Ext$ groups to $\Ext$ groups in every degree, and it preserves projective
covers and the block decomposition because these are defined by properties
invariant under equivalence.  Finally $\Mat_N(B)e\otimes_BV$ has dimension
$N\dim V$, so composing $r$ steps gives $\dim L=p^r\dim L_r$.
\end{proof}

\section{A Peirce--Tietze presentation theorem}\label{sec:peirce-tietze}

\begin{theorem}[Peirce--Tietze corner presentation]
\label{thm:PT}
Let
\[
 A=k\langle x_1,\ldots,x_m\rangle/(R_1(x),\ldots,R_q(x))
\]
be a finitely presented unital associative $k$-algebra, and suppose
noncommutative words $w_{ab}(x)$, $1\le a,b\le N$, have images $E_{ab}\in A$
satisfying
\[
 E_{ab}E_{cd}=\delta_{bc}E_{ad},\qquad\sum_{a=1}^NE_{aa}=1_A,
\]
so that the $E_{ab}$ form a same-unit copy of $\Mat_N(k)$.  Put $e=E_{11}$
and $C=eAe$.  Introduce free noncommuting variables $z_{s;ab}$
($1\le s\le m$, $1\le a,b\le N$), write $F=k\langle z_{s;ab}\rangle$, set
$Z_s=(z_{s;ab})_{a,b}\in\Mat_N(F)$, let $\epsilon_{ab}$ denote the constant
matrix units in $\Mat_N(F)$, and let $I_{\PT}\triangleleft F$ be the ideal
generated by every matrix entry of
\[
 \begin{aligned}
  R_\lambda(Z_1,\ldots,Z_m)&=0 &&(1\le\lambda\le q),\\
  w_{ab}(Z_1,\ldots,Z_m)&=\epsilon_{ab} &&(1\le a,b\le N).
 \end{aligned}
\]
Put $B_{\PT}=F/I_{\PT}$.  Then the assignment
\[
 z_{s;ab}\longmapsto E_{1a}x_sE_{b1}
\]
induces a chart-relative isomorphism $\psi:B_{\PT}\xrightarrow{\sim}eAe$.
Consequently the displayed Peirce coefficients generate $eAe$, and the
relator and chart equations above generate the full ideal of relations
among them.
\end{theorem}

\begin{proof}
The relator equations define a homomorphism
$\Phi:A\to\Mat_N(B_{\PT})$, $x_s\mapsto Z_s$.  The chart equations give
$\Phi(E_{ab})=\epsilon_{ab}$, so $\epsilon_{1a}Z_s\epsilon_{b1}=
z_{s;ab}\epsilon_{11}\in\operatorname{im}\Phi$.  Products show
$b\epsilon_{11}\in\operatorname{im}\Phi$ for every $b\in B_{\PT}$;
multiplying by the constant matrix units gives $b\epsilon_{rs}\in
\operatorname{im}\Phi$ for all $r,s$, so $\Phi$ is surjective.

For $a\in A$ define $\Theta(a)=(E_{1r}aE_{s1})_{r,s}\in\Mat_N(C)$.
Insertion of $1=\sum_tE_{tt}$ proves multiplicativity, and
$(c_{rs})\mapsto\sum_{r,s}E_{r1}c_{rs}E_{1s}$ is its inverse, so $\Theta$
is an algebra isomorphism with $\Theta(E_{ab})=\epsilon_{ab}$.  The
matrices $\Theta(x_s)$ satisfy the relator and chart equations, so
$z_{s;ab}\mapsto E_{1a}x_sE_{b1}$ defines $\psi:B_{\PT}\to C$, and direct
evaluation on every $x_s$ gives
\[
 \Mat_N(\psi)\circ\Phi=\Theta.
\]
Since $\Theta$ is injective, this makes $\Phi$ injective, hence an
isomorphism; then $\Mat_N(\psi)=\Theta\Phi^{-1}$ is an isomorphism, and
taking the $(1,1)$ corner shows $\psi$ is an isomorphism.  
\end{proof}

\begin{remark}
The presentation is relative to the fixed word chart $\{w_{ab}\}$.
Conjugating the chart transports its Peirce generators and relations
along the induced corner isomorphism.
\end{remark}

For a reduced enveloping algebra, the defining relations are indexed
by a fixed Lie algebra basis:

\begin{lemma}[Basis $p$-relations generate the reduced ideal]
\label{lem:p-relations}
Let $\mathfrak g$ be a finite-dimensional restricted Lie algebra, let
$X_1,\ldots,X_d$ be a basis, and write
$[X_i,X_j]=\sum_rc_{ij}^rX_r$, $X_i^{[p]}=\sum_rd_i^rX_r$.  Then
$u_\chi(\mathfrak g)$ is presented by the generators $X_i$, the Lie
relators $X_iX_j-X_jX_i-\sum_rc_{ij}^rX_r$, and the $d$ basis relators
$X_i^p-\sum_rd_i^rX_r-\chi(X_i)^p1$.
\end{lemma}

\begin{proof}
After imposing the Lie relators one obtains $U(\mathfrak g)$.  The
restricted-envelope map
$\xi(x)=x^p-x^{[p]}\in Z(U(\mathfrak g))$ is Frobenius-semilinear:
\[
 \xi\!\left(\sum_i a_iX_i\right)=\sum_i a_i^p\xi(X_i),
 \qquad
 \chi\!\left(\sum_i a_iX_i\right)^p=\sum_i a_i^p\chi(X_i)^p.
\]
Thus the reduced relation for every $x\in\mathfrak g$ is the same
Frobenius-linear combination of the displayed basis relations, and the two
defining ideals coincide.
\end{proof}

\begin{corollary}[Reduced corner presentation]\label{cor:reduced-corner}
If explicit words in the basis generators of $u_\chi(\mathfrak g)$ form a
same-unit $\Mat_N(k)$, the primitive corner has the presentation of
\cref{thm:PT}, with the matrix coefficients of the Lie relators, the basis
$p$-relators, and the word-chart equations.
\end{corollary}

\begin{proof}
The hypothesis places a same-unit $\Mat_N(k)$ inside
$u_\chi(\mathfrak g)$, so $u_\chi(\mathfrak g)\cong\Mat_N(eu_\chi(\mathfrak
g)e)$ for a diagonal matrix unit $e$, and
\cref{prop:morita-transport} identifies the two module categories.  Applying
\cref{thm:PT} to this presentation of $u_\chi(\mathfrak g)$ returns a
presentation of the corner $eu_\chi(\mathfrak g)e$ whose relators are the
images of the defining relators, that is the matrix coefficients of the Lie
relators and of the basis $p$-relators together with the equations
expressing the chosen words as matrix units.
\end{proof}

\section{Iterated Peirce coefficients}\label{sec:iterated}

By \cref{cor:reduced-corner}, every corner produced by the descent of
\cref{thm:uniform-descent} is presented by
the matrix coefficients of a fixed generating set of $\mathfrak g$ against the
chosen chart at each stage.

Fix a sequence of unital matrix reductions
\[
 A_0=u_\chi(\mathfrak g),
 \qquad A_{j+1}=e_jA_je_j,
 \qquad
 M_j=\spn\{E^{(j)}_{ab}\}\cong\Mat_p(k)\subseteq A_j.
\]

\begin{theorem}[Iterated matrix units and Peirce coefficients]\label{thm:composite-coeff}
There are composite matrix units
\[
 \mathbb E^{[r]}_{\mathbf i,\mathbf j}\in A_0,
 \qquad \mathbf i,\mathbf j\in\{0,\ldots,p-1\}^r,
\]
whose span
\[
 M^{[r]}:=\spn\{\mathbb E^{[r]}_{\mathbf i,\mathbf j}\}
 \cong\Mat_{p^r}(k)
\]
is unital in $A_0$, with
\[
 e^{[r]}=\mathbb E^{[r]}_{\mathbf 0,\mathbf 0},
 \qquad A_r=e^{[r]}A_0e^{[r]}.
\]
The coefficient maps
\[
 \gamma^{[r]}_{\mathbf i,\mathbf j}(x)
 =\mathbb E^{[r]}_{\mathbf0,\mathbf i}
 x\mathbb E^{[r]}_{\mathbf j,\mathbf0}
\]
satisfy the composition rule
\[
 \gamma^{[r+1]}_{(\mathbf i,a),(\mathbf j,b)}
 =\gamma^{(r)}_{ab}\circ\gamma^{[r]}_{\mathbf i,\mathbf j}
\]
and reconstruction
\[
 x=\sum_{\mathbf i,\mathbf j}
 \mathbb E^{[r]}_{\mathbf i,\mathbf0}
 \gamma^{[r]}_{\mathbf i,\mathbf j}(x)
 \mathbb E^{[r]}_{\mathbf0,\mathbf j}.
\]
If $G$ is a finite associative generating set of $A_0$, then
\[
 A_r=k\Big\langle
 \gamma^{[r]}_{\mathbf i,\mathbf j}(g):
 g\in G
 \Big\rangle.
\]
\end{theorem}

\begin{proof}
Define recursively
\[
 \mathbb E^{[1]}_{i,j}=E^{(0)}_{ij},
\]
and
\[
 \mathbb E^{[r+1]}_{(\mathbf i,a),(\mathbf j,b)}
 =\mathbb E^{[r]}_{\mathbf i,\mathbf0}
 E^{(r)}_{ab}
 \mathbb E^{[r]}_{\mathbf0,\mathbf j}.
\]
Since $E^{(r)}_{ab}$ lies in the $(\mathbf0,\mathbf0)$ corner, direct
multiplication gives the matrix-unit relations and the displayed coefficient
recursion.  The reconstruction formula follows by inserting the identities
of the successive matrix factors.  Finally, coefficients of products are
polynomials in coefficients of the factors, so coefficients of words in $G$
are generated by coefficients of $G$ itself.  An element of $A_r$ is its own
$(\mathbf0,\mathbf0)$ composite coefficient.
\end{proof}

A finite Lie generating set of $\mathfrak g$ generates $u_\chi(\mathfrak g)$
as an associative algebra, so its Peirce coefficients generate every corner; for $\mathfrak g=\kk$ this is the Contact
algebra of \cref{sec:contact-setup}, and for $\mathfrak
g\in\{\mathfrak w,\mathfrak s,\mathfrak h\}$ it is one of the families of
\cref{sec:wsh}.

\begin{proposition}[Lie return from a Peirce corner]\label{prop:lie-return}
Let $\mathfrak g$ be any restricted simple Cartan-type Lie algebra with
$A_0=u_\chi(\mathfrak g)$.  Let $I\triangleleft A_r$ and let
\[
 \widetilde I^{[r]}
 =\sum_{\mathbf i,\mathbf j}
 \mathbb E^{[r]}_{\mathbf i,\mathbf0}
 I
 \mathbb E^{[r]}_{\mathbf0,\mathbf j}
 \triangleleft A_0
\]
be its ambient Morita lift.  Suppose $\chi\ne0$ and there exists
$0\ne x\in\mathfrak g\cap C_{A_0}(M^{[r]})$ such that
$e^{[r]}xe^{[r]}\in I$.  Then $I=A_r$.
\end{proposition}

\begin{proof}
Since $x$ centralizes $M^{[r]}$, the matrix-unit identities give
\[
 x=\sum_{\mathbf i}
 \mathbb E^{[r]}_{\mathbf i,\mathbf0}
 (e^{[r]}xe^{[r]})
 \mathbb E^{[r]}_{\mathbf0,\mathbf i}.
\]
Thus $x\in\widetilde I^{[r]}$.
Hence $\widetilde I^{[r]}\cap\mathfrak g$ is a nonzero Lie ideal of the
simple Lie algebra $\mathfrak g$, so it equals $\mathfrak g$.  Since $\chi\ne0$, choose
$y\in\mathfrak g$ with $\chi(y)\ne0$.  The reduced relation
\[
 y^p-y^{[p]}=\chi(y)^p1
\]
shows that the two-sided ideal generated by $\mathfrak g$ contains $1$.
Therefore $\widetilde I^{[r]}=A_0$, and taking the composite corner gives
$I=A_r$.
\end{proof}

\section{The intrinsic \texorpdfstring{$p$}{p}-Peirce index}\label{sec:index}

\begin{definition}
For a finite-dimensional $k$-algebra $C$ with
\[
 C/\Jac(C)\cong\prod_{\lambda=1}^t\Mat_{n_\lambda}(k),
\]
define
\[
 \nuP(C)=\min_\lambda\vp(n_\lambda).
\]
Call $C$ \emph{$p$-primitive} if $\nuP(C)=0$.
\end{definition}

\begin{proposition}[Maximal unital matrix factor]\label{prop:max-factor}
A finite-dimensional algebra $C$ contains a unital
$\Mat_{p^s}(k)$ if and only if $s\le\nuP(C)$.  Consequently
\[
 C\cong\Mat_{p^{\nuP(C)}}(D)
\]
for a $p$-primitive full corner $D$.
\end{proposition}

\begin{proof}
If $\Mat_{p^s}(k)$ embeds with the same unit, every simple $C$-module
restricts to a direct sum of the unique $p^s$-dimensional simple module for
that matrix algebra, so $p^s\mid n_\lambda$ for every $\lambda$.
Conversely, Wedderburn--Malcev gives a semisimple complement
$S\cong C/\Jac(C)$.  If $p^s\mid n_\lambda$ for all $\lambda$, diagonally
embed $\Mat_{p^s}(k)$ into every factor by
$X\mapsto X\otimes I_{n_\lambda/p^s}$.  A primitive matrix corner then gives
the asserted matrix decomposition and divides every simple dimension by
$p^s$.
\end{proof}

\begin{proposition}[Peirce index under matrix factors]\label{prop:index-drop}
If $C\cong\Mat_{p^r}(D)$, then
\[
 \nuP(C)=r+\nuP(D).
\]
In particular, each unital $\Mat_p$ reduction lowers $\nuP$ by one.
\end{proposition}

\begin{proof}
The semisimple factor sizes of $C$ are $p^r$ times those of $D$.
\end{proof}

\begin{lemma}[Commutator-unit criterion]\label{lem:comm-unit}
For
\[
 S=\prod_\lambda\Mat_{n_\lambda}(k)
\]
one has
\[
 1_S\in[S,S]
 \iff p\mid n_\lambda\ \text{for every }\lambda.
\]
Consequently, for any finite-dimensional algebra $C$ with radical $J$,
\[
 C\text{ contains a unital }\Mat_p(k)
 \iff 1_C\in J+[C,C].
\]
\end{lemma}

\begin{proof}
For a full matrix algebra,
\[
 [\Mat_n(k),\Mat_n(k)]=\ker\operatorname{tr}.
\]
Indeed off-diagonal matrix units and diagonal differences are commutators.
Thus $I_n$ is a commutator exactly when $n=0$ in $k$, that is, $p\mid n$.
Passing to the semisimple quotient gives the second assertion together with
\cref{prop:max-factor}.
\end{proof}

\begin{theorem}[Central-socle trace space]\label{thm:trace-space}
Let $C$ be a finite-dimensional symmetric algebra with symmetrizing
functional $\tau$, and put $J=\Jac(C)$, $S=C/J$ and
$R=Z(C)\cap\Soc(C)$. For $s\in R$ define
\[
 \ell_s:S\to k,
 \qquad \ell_s(\bar x)=\tau(xs).
\]
Then
\[
 R\xrightarrow{\sim}\Tr(S)
 :=\{\ell\in S^*: \ell(xy)=\ell(yx)\}
 =(S/[S,S])^*.
\]
There exists $s_0\in R$ such that
\[
 \beta_{s_0}(\bar x,\bar y)=\tau(xys_0)
\]
is nondegenerate on $S$.  For this $s_0$, every $s\in R$ corresponds to a
unique $z_s\in Z(S)$ satisfying
\[
 \ell_s(x)=\ell_{s_0}(z_sx),
\]
and
\[
 R\xrightarrow{\sim}Z(S),\qquad s\mapsto z_s.
\]
\end{theorem}

\begin{proof}
For the nondegenerate pairing $\langle x,y\rangle=\tau(xy)$,
$J^\perp=\Soc(C)$ and $[C,C]^\perp=Z(C)$.
Indeed, orthogonality to $J$ is equivalent to $Jx=0$;
orthogonality to all commutators gives
$\tau(a(bx-xb))=0$ for every $a,b$, hence $bx=xb$.
Consequently
\[
 R=(J+[C,C])^\perp.
\]
The functional $\ell_s$ is well defined because $Js=0$, and is a trace
because $s$ is central. The induced pairing
\[
 S/[S,S]\times R\to k,
 \qquad ([\bar x],s)\mapsto\tau(xs)
\]
is perfect; hence $R\cong\Tr(S)$.

A semisimple algebra admits a nondegenerate trace form: after an
Artin--Wedderburn decomposition, take the sum of the ordinary matrix traces.
The form $\operatorname{tr}(xy)$ on $\Mat_n(k)$ is nondegenerate.
Its inverse image determines a unique $s_0\in R$.  For arbitrary $s$,
nondegeneracy gives a unique $z_s$
with the displayed identity.  For $a,x\in S$,
\[
 \ell_{s_0}((z_sa-az_s)x)
 =\ell_s(ax)-\ell_s(xa)=0,
\]
so $z_s$ is central.  Conversely multiplication of $\ell_{s_0}$ by any
central $z$ is again a trace functional, hence arises from a unique $s$.
\end{proof}

\begin{theorem}[Central-socle detector]\label{thm:central-socle-detector}
Let $C$ be a finite-dimensional symmetric algebra with symmetrizing
functional $\tau$.  Then
\begin{align*}
 C\text{ contains a unital }\Mat_p(k)
 &\Longleftrightarrow
 \tau|_{Z(C)\cap\Soc(C)}=0.
\end{align*}
Equivalently, $C$ is $p$-primitive if and only if
$\tau$ is nonzero on the central socle.
\end{theorem}

\begin{proof}
By \cref{thm:trace-space}, every trace on $C/J$ is $\ell_s$ for a
unique $s\in R$, and $\tau(s)=\ell_s(1)$.
The ordinary trace on the factor $\Mat_{n_\lambda}(k)$ takes value
$n_\lambda$ at its identity. Thus $\tau$ vanishes on $R$ precisely
when $p$ divides every $n_\lambda$. Apply \cref{prop:max-factor}.
\end{proof}

\begin{proof}[Proof of \cref{thm:uniform-descent}]
A full corner of a symmetric algebra is symmetric.
By \cref{thm:central-socle-detector}, a unital $\Mat_p(k)$ exists
exactly when the symmetrizing functional vanishes on the central socle.
Each reduction $C_j\cong\Mat_p(C_{j+1})$ lowers $\nuP$ by one,
by \cref{prop:index-drop}. The process ends after $\nuP(C_0)$ steps,
and \cref{prop:morita-transport} gives the module equivalences.
\end{proof}

\section{Central-socle trace pencils}\label{sec:trace-pencil}

Let $C$ be a finite-dimensional symmetric algebra, with
$J=\Jac(C)$, $S=C/J$ and $R=Z(C)\cap\Soc(C)$.
Fix $s_0$ and the identification $R\cong Z(S)$ of \cref{thm:trace-space}.

\begin{theorem}[One generic trace separates all factors]\label{thm:generic-trace}
Let $t=\dim_kR$.  Choose $s_1\in R$ so that
$z_1=z_{s_1}\in Z(S)$ has squarefree minimal polynomial of degree $t$,
\[
 \mu_{z_1}(X)=\prod_{\alpha=1}^t(X-\lambda_\alpha),
 \qquad \lambda_\alpha\ne\lambda_\beta\quad(\alpha\ne\beta).
\]
Then
\[
 z_\alpha
 =\prod_{\beta\ne\alpha}
 \frac{z_1-\lambda_\beta}{\lambda_\alpha-\lambda_\beta}
\]
are precisely the primitive central idempotents of $S$, and
\[
 S=\prod_{\alpha=1}^tS_\alpha,
 \qquad
 S_\alpha=z_\alpha S\cong\Mat_{n_\alpha}(k),
 \qquad
 n_\alpha=\sqrt{\dim_kS_\alpha}.
\]
Equivalently, for
\[
 \beta_\lambda(x,y)=\ell_{s_1}(xy)-\lambda\ell_{s_0}(xy),
\]
one has
\[
 \rad\beta_{\lambda_\alpha}=S_\alpha.
\]
\end{theorem}

\begin{proof}
Since $Z(S)$ is a finite product of copies of the infinite field $k$, it
contains an element with pairwise distinct coordinates.
The Lagrange polynomials are then exactly the coordinate projectors.  On the
$\alpha$-factor, $\beta_{\lambda_\alpha}$ vanishes, while on every other
factor it is a nonzero scalar multiple of a nondegenerate matrix trace form.
Thus its radical is precisely $S_\alpha$.
\end{proof}

\section{Trace reconstruction and extension bimodules}\label{sec:extensions}

Let $C$ be the final $p$-primitive corner, and let
\[
 G_C=\{g_1,\ldots,g_N\}\subset C
\]
be the finite generating set provided by \cref{thm:composite-coeff}.

\begin{proposition}[The base trace form recovers the radical]\label{prop:recover-radical}
For the $s_0$ of \cref{thm:trace-space}, the form
\[
 B_{s_0}(x,y)=\tau(xys_0)
\]
on $C$ satisfies
\[
 \rad B_{s_0}=\Jac(C).
\]
\end{proposition}

\begin{proof}
The inclusion $J\subseteq\rad B_{s_0}$ follows from $Js_0=0$.  If
$x\notin J$, then $\bar x\ne0$ in $S$; nondegeneracy of the descended form
on $S$ gives $\bar y$ with $\ell_{s_0}(\bar x\bar y)\ne0$.  Any lift $y$
then satisfies $B_{s_0}(x,y)\ne0$.
\end{proof}

\begin{proposition}[Generator actions on the simple factors]\label{prop:factor-actions}
Let $z_\alpha$ and $S_\alpha$ be the factors obtained in
\cref{thm:generic-trace}.  Define
\[
 \rho_\alpha(g)=z_\alpha(g+J)\in S_\alpha,
 \qquad g\in G_C.
\]
Then
\[
 S_\alpha=k\langle\rho_\alpha(g):g\in G_C\rangle_{z_\alpha}.
\]
The tuple
\[
 \bigl(z_\alpha,n_\alpha,\rho_\alpha|_{G_C}\bigr)
\]
records the action of the chosen Peirce generators on the unique
simple left $S_\alpha$-module $L_\alpha$.
\end{proposition}

\begin{proof}
$G_C$ generates $C$, hence its image generates $S=C/J$.
Multiplication by the central idempotent $z_\alpha$ projects every word in
these images to the same word in the projected generators, with identity
$z_\alpha$.
\end{proof}

\begin{theorem}[The first radical layer and $\Ext^1$]\label{thm:radical-ext}
Let $M=J/J^2$ and
\[
 M_{\beta\alpha}=z_\beta Mz_\alpha.
\]
As an $(S_\beta,S_\alpha)$-bimodule, $M_{\beta\alpha}$ is a direct sum of
copies of
\[
 B_{\beta\alpha}=L_\beta\otimes_kL_\alpha^*.
\]
The multiplicity space
\[
 \mathsf E_{\alpha\to\beta}
 =\Hom_{S_\beta-S_\alpha}
 (B_{\beta\alpha},M_{\beta\alpha})
\]
satisfies
\[
 M_{\beta\alpha}
 \cong B_{\beta\alpha}\otimes_k\mathsf E_{\alpha\to\beta}
\]
and
\[
 \dim_k\Ext_C^1(L_\alpha,L_\beta)
 =\dim_k\mathsf E_{\alpha\to\beta}.
\]
The blocks are the connected components of the underlying undirected graph
of the directed graph having an arrow $\alpha\to\beta$ whenever
$\mathsf E_{\alpha\to\beta}\ne0$.
\end{theorem}

\begin{proof}
Since
\[
 S_\beta\otimes_kS_\alpha^{\rm op}
 \cong\Mat_{n_\beta n_\alpha}(k),
\]
it has one simple module, namely $L_\beta\otimes L_\alpha^*$.
Semisimplicity gives the displayed decomposition.  Choosing primitive
idempotents in the two factors identifies the multiplicity space with
\[
 e_\beta(J/J^2)e_\alpha,
\]
the primitive-idempotent realization of
$\Ext_C^1(L_\alpha,L_\beta)$.
\end{proof}

\begin{proposition}[Higher extensions]\label{prop:higher-ext}
Choose a Wedderburn--Malcev lift $E\subset C$ of the recovered semisimple
quotient $S$.  For a simple $L_\alpha$, the normalized relative bar terms
\[
 C\otimes_EJ^{\otimes_Eq}\otimes_EL_\alpha,
 \qquad q\ge0,
\]
form a projective resolution computing all higher $\Ext_C^q$.
\end{proposition}

\begin{proof}
The semisimple algebra $E$ is separable over $k$, so the normalized relative
bar complex is a projective resolution.  Its differential is induced by the
multiplication on $J$.
\end{proof}

\begin{proof}[Proof of \cref{thm:trace-reconstruction}]
Combine \cref{thm:trace-space,thm:generic-trace,prop:recover-radical,prop:factor-actions,thm:radical-ext,prop:higher-ext}.
\end{proof}

\section{Descent and reconstruction}\label{sec:main}

\begin{theorem}[$p$-Peirce descent]
\label{thm:main}\label{thm:global}
Let $k$ be algebraically closed of characteristic $p\ge5$, let
\[
 \mathfrak g\in\{W(n),\ S(n)^{(1)}\ (n\ge3),\
 H(2m)^{(2)},\ K(2m+1;\mathbf1)^{(1)}\ (m\ge1)\},
\]
and let $0\ne\chi\in\mathfrak g^*$.
There is a $p$-primitive symmetric algebra $C_\chi$ such that
\[
 u_\chi(\mathfrak g)\cong\Mat_{p^r}(C_\chi),\qquad
 r=\nuP(u_\chi(\mathfrak g))\ge1.
\]
The central-socle trace pencil recovers $J=\Jac(C_\chi)$ and
$S=C_\chi/J=\prod_\alpha S_\alpha$, where
$S_\alpha\cong\Mat_{n_\alpha}(k)$ with identity $z_\alpha$.
The Peirce generators determine the action on each simple module
$L_\alpha$, and
\[
 \dim\Ext^1_{C_\chi}(L_\alpha,L_\beta)
 =\dim\Hom_{S_\beta-S_\alpha}
 \bigl(L_\beta\otimes L_\alpha^*,\ z_\beta(J/J^2)z_\alpha\bigr).
\]
The relative bar resolution gives all higher extension groups.
These data transport to $u_\chi(\mathfrak g)$ under Morita equivalence.
\end{theorem}

\begin{proof}
The first matrix factor is supplied by
\cref{thm:first-matrix,thm:witt-first-matrix,thm:special-first-matrix,thm:hamiltonian-first-matrix}.
Symmetry follows from \cref{prop:symmetric-general}.
Apply \cref{thm:uniform-descent} and
\cref{thm:trace-reconstruction}, using the generators and presentations
of \cref{thm:PT,thm:composite-coeff}.
Finally apply \cref{prop:morita-transport}.
\end{proof}

\begin{remark}[Choice of matrix units]
The integer $r$ and the Morita class of $C_\chi$ are independent of
the matrix charts. Changing the generic trace element permutes the
primitive central factors.
\end{remark}

\section{Basic algebras and transport}\label{sec:harvest}

Throughout this section, $\mathfrak g$ and $\chi$ are as in \cref{thm:global},
and $C_\chi$ is the $p$-primitive corner produced there, with
$J=\Jac(C_\chi)$ and simple modules $\{L_\alpha\}$.

\subsection{A bound-quiver presentation}

\begin{lemma}[Idempotent lifting]\label{lem:idempotent-lift}
Every primitive idempotent of $C_\chi/J$ lifts to a primitive idempotent of
$C_\chi$, and a complete orthogonal family of primitive idempotents lifts to
one in $C_\chi$.
\end{lemma}

\begin{proof}
Since $C_\chi$ is finite-dimensional, $J^\ell=0$ for some $\ell$, and it
suffices to lift through each square-zero layer $J^r/J^{r+1}$.  If $a$ lifts
an idempotent modulo a square-zero ideal $I$ and $d=a^2-a\in I$, then $d$
commutes with $a$, $d^2=0$, and
\[
 e=a+(1-2a)d
\]
satisfies $e^2=a^2+2a(1-2a)d=(a+d)-2ad=a+(1-2a)d=e$.  Iterating over the
finitely many layers lifts one idempotent; lifting the next inside the
complementary corner $(1-e)C_\chi(1-e)$ preserves orthogonality.  A lift of
a primitive idempotent stays primitive, since a nontrivial orthogonal
decomposition would descend to $C_\chi/J$, while $J$ contains no nonzero
idempotent.
\end{proof}

Fix one lifted primitive idempotent $e_\alpha$ above each simple
$L_\alpha$, put $e=\sum_\alpha e_\alpha$, and let
$B_\chi=eC_\chi e$. Since $e$ meets every simple module, it is full;
thus $B_\chi$ is basic and Morita equivalent to $C_\chi$ and
$u_\chi(\mathfrak g)$.  Write
$J_B=\Jac(B_\chi)=eJe$ and let $Q_\chi$ be the Gabriel quiver with one
vertex per $L_\alpha$ and, for each ordered pair $(\alpha,\beta)$, one
arrow for each vector in a chosen basis of
$e_\beta J_Be_\alpha/e_\beta J_B^2e_\alpha$. Choose lifts in
$e_\beta J_Be_\alpha$ for these arrows; they define an algebra map
$\Phi_\chi:kQ_\chi\to B_\chi$.

\begin{theorem}[Bound-quiver presentation]\label{thm:bound-quiver}
Let $\ell$ be the Loewy length of $B_\chi$.  Then $\Phi_\chi$ is
surjective, its kernel is generated by the finite-dimensional space
\[
 R_{<\ell}=\ker\bigl(\Phi_\chi|_{P_{<\ell}}\bigr)
\]
of relations among paths of length $<\ell$ together with the space
$P_\ell$ of all paths of length $\ell$, and
\[
 B_\chi\cong kQ_\chi/\langle R_{<\ell},P_\ell\rangle.
\]
For $P_\alpha=B_\chi e_\alpha$,
\[
 \rad^qP_\alpha/\rad^{q+1}P_\alpha=J_B^qe_\alpha/J_B^{q+1}e_\alpha,
 \qquad
 C_{\beta\alpha}=[P_\alpha:L_\beta].
\]
For every fixed $n\ge0$,
$\Ext^n_{u_\chi(\mathfrak g)}(L_\alpha,L_\beta)$ is the degree-$n$
cohomology of a complex of the projectives $B_\chi e_\gamma$.
\end{theorem}

\begin{proof}
Modulo $J_B$, the idempotents $e_\alpha$ span $B_\chi/J_B$, and the chosen
arrow lifts span $J_B/J_B^2$ as a bimodule over it; inductively, if the
image of $\Phi_\chi$ contains representatives modulo $J_B^q$, multiplying
by arrow lifts produces $J_B^q/J_B^{q+1}=(J_B)^q/J_B^{q+1}$, so surjectivity
follows once $q=\ell$.  Every path of length $\ge\ell$ maps to $0$ and lies
in a length-$\ell$ path; hence, for $x\in\ker\Phi_\chi$, splitting
$x=x_{<\ell}+x_{\ge\ell}$ by path length shows $x_{\ge\ell}$ already lies in
the displayed ideal, and $x_{<\ell}$ does by definition of $R_{<\ell}$.
This proves the presentation.  For the Loewy layers, $\rad^qP_\alpha=
J_B^qe_\alpha$ by definition of the radical filtration of a projective
module; the Cartan and $\Ext^n$ statements follow by constructing the initial segment through degree $n$
of a minimal projective resolution of each $L_\alpha$ over $B_\chi$.  Morita
transport along \cref{prop:morita-transport} carries all of this to
$u_\chi(\mathfrak g)$.
\end{proof}

\begin{corollary}[Second Loewy layers]
\label{cor:loewy-cartan}
For the basic algebra $B_\chi$ and its indecomposable projective
$P_\alpha=B_\chi e_\alpha$, one has, for every ordered pair
$(\alpha,\beta)$,
\[
 \left[\frac{\rad P_\alpha}{\rad^2P_\alpha}:L_\beta\right]
 =\dim_k e_\beta(J_B/J_B^2)e_\alpha
 =\dim_k\Ext^1_{u_\chi(\mathfrak g)}(L_\alpha,L_\beta).
\]
The same equalities hold after transport between
$B_\chi$, $C_\chi$, and $u_\chi(\mathfrak g)$, with the simple modules
relabeled by the corresponding Morita equivalences.
\end{corollary}

\begin{proof}
Since $\rad P_\alpha=J_Be_\alpha$, the first radical layer is
$(J_B/J_B^2)e_\alpha$.  Its $L_\beta$-multiplicity is
\[
 \dim_k e_\beta(J_B/J_B^2)e_\alpha,
\]
which is the primitive-idempotent
realization of $\Ext^1(L_\alpha,L_\beta)$ and agrees with
\cref{thm:radical-ext}.  Morita equivalence preserves projective
covers, radical layers, composition multiplicities, and Yoneda extensions.
\end{proof}

\subsection{Transport under restricted automorphisms}

Let $\phi\in\operatorname{Aut}_{[p]}(\mathfrak g)$, and set
$(\phi\cdot\chi)(y)=\chi(\phi^{-1}y)$ and $\psi=\phi\cdot\chi$.

\begin{proposition}[Automorphism transport]\label{prop:orbit-compression}
The universal-enveloping automorphism $U(\phi)$ carries the reduced ideal
$I_\chi$ onto $I_\psi$, inducing an algebra isomorphism
$\Phi_\phi:u_\chi(\mathfrak g)\xrightarrow{\sim}u_\psi(\mathfrak g)$.
Under $\Phi_\phi$, the data of \cref{thm:global,thm:bound-quiver} for
$\chi$ transport to those for $\psi$: simple modules correspond, the basic
algebra and Gabriel quiver identify after relabelling vertices, every
Loewy layer and the Cartan matrix agree, and
$\Ext^n_{u_\chi(\mathfrak g)}(L,M)\cong
\Ext^n_{u_\psi(\mathfrak g)}(\Phi_{\phi,*}L,\Phi_{\phi,*}M)$ for every
$n\ge0$.  In particular, every invariant of \cref{thm:global,thm:bound-quiver}
is constant, up to this transport, on $\operatorname{Aut}_{[p]}(\mathfrak
g)$-orbits in $\mathfrak g^*$.
\end{proposition}

\begin{proof}
Since $\phi$ is restricted, $\phi(x^{[p]})=\phi(x)^{[p]}$ for $x\in
\mathfrak g$, so with $r_\chi(x)=x^p-x^{[p]}-\chi(x)^p1$,
\[
 U(\phi)(r_\chi(x))
 =\phi(x)^p-\phi(x)^{[p]}-\chi(x)^p1
 =r_\psi(\phi(x)),
\]
using $\chi(x)=\psi(\phi(x))$.  Hence $U(\phi)(I_\chi)\subseteq I_\psi$,
and the reverse inclusion follows by applying the same argument to
$\phi^{-1}$; this gives $\Phi_\phi$.  An algebra isomorphism induces an
exact equivalence of module categories preserving simples, projective
covers, radical and Loewy filtrations, and Yoneda extension groups, which
gives the transport of $\Ext^n$.  Applying $\Phi_\phi$ to the primitive
idempotents $e_\alpha$ of \cref{thm:bound-quiver} gives primitive
idempotents representing the transported simples, and multiplicativity of
$\Phi_\phi$ carries the path map $\Phi_\chi$ to the corresponding path map
for $\psi$ on the relabelled quiver, transporting the presentation of
\cref{thm:bound-quiver} exactly.
\end{proof}

\section{The Witt algebra \texorpdfstring{$W(1;\mathbf1)$}{W(1;1)}}\label{sec:w11}

\subsection{Matrix depth and terminal corners}\label{subsec:w11}

\begin{example}\label{ex:w11}
Let $\mathfrak w=W(1;\mathbf1)$ with basis $e_i=x^{i+1}\partial$,
$-1\le i\le p-2$.  Then $[e_i,e_j]=(j-i)e_{i+j}$, $e_0^{[p]}=e_0$ and
$e_i^{[p]}=0$ for $i\ne0$; here $r=p-2$ and $\mathfrak w_{(j)}$ is spanned
by the $e_i$ with $i\ge j$.  Let $0\ne\chi\in\mathfrak w^*$ and put
$h=h(\chi)$.

\emph{First matrix factor.}  The proof of \cref{thm:witt-first-matrix}
gives the pair.  If $a=\chi(e_{-1})\ne0$, take $U=a^{-1}e_{-1}$
and $h_0=-e_0$; then $U^p=1$ in $u_\chi(\mathfrak w)$ and $[h_0,U]=U$.  If
$\chi(e_{-1})=0$ and $s\ge1$ is maximal with $\chi(e_s)\ne0$, take
$\Theta=e_s$ and $h_0=s^{-1}e_0$, so that $[h_0,\Theta]=\Theta$; here
$1\le s\le p-2$, so $s$ is invertible in $k$.  If $\chi$ is supported on
$e_0$ alone, the substitution $x\mapsto x+cx^2$ with $c\ne0$ moves $\chi$
to a character with $\chi(e_{-1})\ne0$.

\emph{Depth and terminal corner.}  By \cref{thm:uniform-descent},
$u_\chi(\mathfrak w)\cong\Mat_{p^\nu}(C_\chi)$, where $\nu$ is the minimum
of the $p$-adic valuations of the dimensions of the simple modules and
$C_\chi$ is symmetric with no unital copy of $\Mat_p(k)$.  The simple
$u_\chi(\mathfrak w)$-modules, and their dimensions, are given by Chang
\cite{Chang1941}, by Feldvoss and Nakano \cite{FeldvossNakano1998}, and
at nonzero height by Boe, Nakano and Wiesner
\cite[\S4]{BoeNakanoWiesner2009}: for $h=0$ there are $p-1$ of them and for
$h=1$ there are $p$, all of dimension $p$; for $1<h<p-1$ there is a unique
simple module $L$, of dimension $p^{\lfloor h/2\rfloor+1}$; and for $h=p-1$
all simple modules have dimension $p^{(p-1)/2}$, with $p$ of them when
$\mathfrak w^\chi$ is a torus and $p-1$ of them when $\mathfrak w^\chi$ is
$p$-nilpotent.  Since $\dim u_\chi(\mathfrak w)=p^p$, this gives:
\begin{enumerate}
\item for $h\in\{0,1\}$, $\nu=1$ and $\dim C_\chi=p^{p-2}$;
\item for $1<h<p-1$, $\nu=\lfloor h/2\rfloor+1$ and $C_\chi$ is a local
symmetric algebra of dimension $p^{\,p-2\lfloor h/2\rfloor-2}$;
\item for $h=p-1$, $\nu=(p-1)/2$ and $\dim C_\chi=p$; then
$C_\chi\cong k^p$ in the toral case, and
$C_\chi\cong k^{p-2}\times k[t]/(t^2)$ in the $p$-nilpotent case, by
Feldvoss and Nakano \cite[Theorem~2.6]{FeldvossNakano1998}.
\end{enumerate}

\emph{Local corners.}  For $1<h<p-1$ the corner $C_\chi$ is local,
and \cref{thm:trace-reconstruction} gives
\[
 \dim\Ext^1_{u_\chi(\mathfrak w)}(L,L)
 =\dim_k\Jac(C_\chi)/\Jac(C_\chi)^2 ,
\]
the embedding dimension of $C_\chi$.  In particular
$\Ext^1_{u_\chi(\mathfrak w)}(L,L)\cong k$ exactly when
$C_\chi\cong k[t]/(t^{\dim C_\chi})$.  For $p=5$ and $h\in\{2,3\}$ one has
$\dim C_\chi=5$, and Boe, Nakano and Wiesner showed that $u_\chi(\mathfrak w)$
is uniserial \cite[\S4.3]{BoeNakanoWiesner2009}; hence
$C_\chi\cong k[t]/(t^5)$. For even $h$ with $2\le h\le p-3$, \cref{thm:bnw-main} below gives
\[
 \dim C_\chi=p^{p-h-2},\qquad
 \dim\Jac(C_\chi)/\Jac(C_\chi)^2=\min\{p-h-2,h+2\}.
\]
The Ext quiver has one vertex and $\min\{p-h-2,h+2\}$ loops.
Its embedding dimension is one precisely when $h=p-3$; in that case
$C_\chi\cong k[t]/(t^p)$. A finite-dimensional local algebra
with residue field $k$ and one-dimensional $J/J^2$ is generated by
one nilpotent element, so it is $k[t]/(t^N)$ with $N=\dim C_\chi$.
Thus $h=p-3$ is the only even height in this range with a uniserial block.
For $p\ge2h+5$ the embedding dimension is $h+2$, and
\[
 \dim\Jac(C_\chi)^2=p^{p-h-2}-h-3.
\]
\end{example}

\subsection{Self-extensions at even height}\label{subsec:w11-even}

Boe, Nakano and Wiesner conjectured that their lower bound for
self-extensions is attained at every even height
\cite[\S4.3]{BoeNakanoWiesner2009}. \Cref{thm:bnw-main} gives the self-extension dimension and hence
the embedding dimension of the local corner in \cref{ex:w11}.

\begin{theorem}\label{thm:bnw-main}
Let $\mathfrak w=W(1;\mathbf1)$, and let $\chi$ have even height
$2\le h\le p-3$. For its unique simple module $L$,
\[
 \dim\Ext^1_{u_\chi(\mathfrak w)}(L,L)=\min\{p-h-2,h+2\}.
\]
The lower bound of Boe, Nakano and Wiesner
\cite[Proposition~4.3(B)]{BoeNakanoWiesner2009}
is attained precisely when $\mathfrak w_{(h)}$ is abelian, that is,
when $p\le2h+1$.
\end{theorem}

Characters of even height $h$ form a single automorphism orbit
\cite[\S4.3]{BoeNakanoWiesner2009}. Choose
$\chi(e_{h-1})=1$ and all other basis values zero. For the remainder of this section,
write $\mathfrak w_j=\mathfrak w_{(j)}$ and put $s=h/2$.
Let $v$ generate the one-dimensional $u_\chi(\mathfrak w_s)$-module
$k_\chi$, and set
\[
 L=u_\chi(\mathfrak w)\otimes_{u_\chi(\mathfrak w_s)}k_\chi,
 \qquad S=u_\chi(\mathfrak w_0)\otimes_{u_\chi(\mathfrak w_s)}k_\chi.
\]
These are the unique simple modules of the respective algebras
\cite[\S4.3]{BoeNakanoWiesner2009}.

\subsubsection{Cohomology of the induced module}

Let $\beta_\chi$ be the alternating form $\beta_\chi(x,y)=\chi([x,y])$
on $\mathfrak{w}_0$. Its matrix on
$\langle e_0,\ldots,e_{h-1}\rangle$ is anti-diagonal: the entry pairing
$e_i$ with $e_{h-1-i}$ equals $h-1-2i$, and all other entries vanish.
Each anti-diagonal entry is nonzero in $k$, because $h-1$ is odd and
$|h-1-2i|<p$. The form is therefore nondegenerate on that span. On the
other hand $[\mathfrak{w}_h,\mathfrak{w}_0]\subseteq\mathfrak{w}_h$ and $\chi(\mathfrak{w}_h)=0$, so
$\mathfrak{w}_h$ lies in the radical. Hence
$\operatorname{rad}(\beta_\chi)=\mathfrak{w}_h$. Together with the fact that
$\mathfrak{w}_h$ is an ideal of $\mathfrak{w}_0$, this supplies the hypotheses of
\cite[Theorem 4.3(C)]{BoeNakanoWiesner2009}.

The ideal $\mathfrak{w}_h$ kills $S$, so $S$ is a simple module for
$Q=u_{\bar\chi}(\mathfrak{w}_0/\mathfrak{w}_h)$. Then $\dim Q=p^h=p^{2s}$ and
$\dim S=p^s$. The representation $Q\to\operatorname{End}_k(S)$ is
surjective by the density theorem, hence bijective by dimension, so $Q$
is semisimple. \begin{proposition}[Induction and the semisimple quotient]
\label{bnw:prop:adjunction}
For every $n\ge0$,
\[
 \Ext^n_{u_\chi(\mathfrak w)}(L,L)
 \cong\Hom_Q\bigl(S,H^n(u_0(\mathfrak w_h),L)\bigr),
\]
and
\begin{equation}\label{bnw:eq:BNW}
 \dim\Ext^n_{u_\chi(\mathfrak w)}(L,L)
 =p^{-s}\dim H^n(u_0(\mathfrak w_h),L).
\end{equation}
\end{proposition}

\begin{proof}
Transitivity gives
$L=u_\chi(\mathfrak w)\otimes_{u_\chi(\mathfrak w_0)}S$.
The PBW theorem makes this induction exact and preserves projectives,
so adjunction holds in every degree.
For the ideal $\mathfrak w_h\triangleleft\mathfrak w_0$, the
Lyndon--Hochschild--Serre spectral sequence has
\[
 E_2^{i,j}=\Ext_Q^i\bigl(S,H^j(u_0(\mathfrak w_h),L)\bigr).
\]
Semisimplicity of $Q$ leaves only the column $i=0$.
Every $Q$-module is a sum of copies of its $p^s$-dimensional simple
module $S$, giving the dimension identity.
\end{proof}

Write
\[
V=\operatorname{span}\{e_0^{a_0}\cdots e_{s-1}^{a_{s-1}}v:0\le a_i<p\},
\qquad F_m=\bigoplus_{a=0}^{m}e_{-1}^aV,
\]
with $F_m=0$ for $m<0$. The relative $\operatorname{PBW}$ theorem gives $L=F_{p-1}$
and $\dim V=p^s$. Also $\mathfrak{w}_hV=0$: commuting an element $e_j$ with $j\ge h$ through a
monomial in $e_0,\ldots,e_{s-1}$ leaves only such elements, all of
which kill $v$.

\begin{lemma}\label{bnw:lem:pivot}
Let $X=e_h$. Then $\ker X^j=F_{j-1}$ for $1\le j\le p$, and
$\operatorname{im}X=\ker X^{p-1}$. For every $i\ge0$ with
$h+i\le p-2$,
\begin{equation}\label{bnw:eq:degree}
e_{h+i}F_m\subseteq F_{m-i-1}.
\end{equation}
The map $X^j$ lowers the highest nonzero $\operatorname{PBW}$ layer by
exactly $j$ whenever that layer has index at least $j$.
\end{lemma}

\begin{proof}
Commuting $e_{h+i}$ through $e_{-1}^a$, the terms with fewer than
$i+1$ commutators still carry some $e_j$ with $j\ge h$ on the right,
and so kill $V$. The first possibly nonzero term is
\begin{equation}\label{bnw:eq:leading}
e_{h+i}e_{-1}^aw\equiv(-1)^{i+1}\binom{a}{i+1}
\frac{(h+i+1)!}{h!}\,e_{-1}^{a-i-1}e_{h-1}w
\pmod{F_{a-i-2}}
\end{equation}
for $a\ge i+1$ and $w\in V$, which gives \eqref{bnw:eq:degree}, the case
$a<i+1$ included. For $i=0$ the leading coefficient is $-a(h+1)$ and
the leading operator is $e_{h-1}$; both are invertible for
$1\le a<p$, since $h+1<p$ and $e_{h-1}^p=\chi(e_{h-1})^p=1$ in
$u_\chi(\mathfrak{w})$. Hence $X^{p-1}:F_{p-1}/F_{p-2}\to V$ is a composite of
invertible leading maps and has rank $\dim V$. As $X^p=0$ and
$\dim L=p\dim V$, every Jordan block of $X$ has length $p$. The
inclusions $F_{j-1}\subseteq\ker X^j$ are then equalities by dimension,
and the image identity follows from the same blocks.
\end{proof}

Every element of $\mathfrak{w}_h$ has zero restricted $p$-th power, since as a
derivation of $k[x]/(x^p)$ it strictly raises degree. A restricted
$1$-cocycle $d:\mathfrak{w}_h\to L$ therefore satisfies
\begin{equation}\label{bnw:eq:restricted}
d([x,y])=x\,d(y)-y\,d(x),\qquad x^{p-1}d(x)=0\qquad(x,y\in\mathfrak{w}_h).
\end{equation}
Applying the second condition to $X$ and using Lemma~\ref{bnw:lem:pivot}
gives $d(X)\in\operatorname{im}X$, so subtracting an inner cocycle
normalises $d(X)=0$. Two such normalisations differ by an inner cocycle $d_u$ with
$Xu=0$, that is with $u\in V$; since $\mathfrak{w}_hV=0$, such a $d_u$ vanishes
identically. Normalised cocycles therefore correspond bijectively to
classes in $H^1(u_0(\mathfrak{w}_h),L)$.

\subsubsection{The abelian range and $p=2h+3$}

If $p\le2h+1$, then $\mathfrak{w}_h$ is abelian. For $j>h$ the bracket
$[X,e_j]$ vanishes, so $Xd(e_j)=e_jd(X)=0$ and $d(e_j)\in\ker X=V$.
Conversely any values in $V$ on $e_{h+1},\ldots,e_{p-2}$ satisfy
\eqref{bnw:eq:restricted}, because $\mathfrak{w}_hV=0$. Hence
\begin{equation}\label{bnw:eq:abelian}
\dim H^1(u_0(\mathfrak{w}_h),L)=(p-h-2)\dim V.
\end{equation}

Let $p=2h+3$. The only nonzero basis bracket in $\mathfrak{w}_h$ is
$[X,e_{h+1}]=e_{2h+1}$. Put $a=d(e_{h+1})$; then $d(e_{2h+1})=Xa$, and
$[X,e_{2h+1}]=0$ forces $a\in\ker X^2=F_1$. For $2\le j\le h$ the
bracket $[X,e_{h+j}]=j\,e_{2h+j}$ vanishes, since $2h+j>p-2$, so
$Xd(e_{h+j})=0$ and $d(e_{h+j})$ is an arbitrary element of $V$.
Conversely all these values lie in $F_1$, and
$e_{h+j}F_1\subseteq F_{-j}=0$ for $j\ge1$ by \eqref{bnw:eq:degree}, so
every bracket other than $[X,e_{h+1}]$ has both action terms zero;
the restricted condition holds because $\mathfrak{w}_hF_1\subseteq V$ and
$\mathfrak{w}_hV=0$. Hence
\begin{equation}\label{bnw:eq:first}
\dim H^1(u_0(\mathfrak{w}_h),L)=(h+1)\dim V.
\end{equation}

\subsubsection{The range $p\ge2h+5$}

For the remaining primes $p\ge2h+5$ write $a_i=d(e_{h+i})$ for
$1\le i\le h$ and normalise $d(X)=0$. The values on all higher basis elements are determined by these,
since $[X,e_j]=(j-h)e_{h+j}$ generates them successively with
$j-h\not\equiv0$.

\begin{proposition}\label{bnw:prop:lower}
If $p\ge2h+5$, then $\dim H^1(u_0(\mathfrak{w}_h),L)\ge(h+2)\dim V$.
\end{proposition}

\begin{proof}
Take $a_1,a_2\in F_1$ and $a_3,\ldots,a_h\in V$ arbitrary, put
\[
d(e_{2h+1})=Xa_1,\qquad d(e_{2h+2})=\tfrac12Xa_2,
\]
and set all other higher values to zero; the two displayed values lie
in $V$. The brackets with $X$ hold by construction: $[X,e_{h+1}]$ and
$[X,e_{h+2}]$ give the two displayed values, $[X,e_{h+j}]$ with
$j\ge3$ gives $Xa_j=0$, and $[X,e_{2h+t}]$ gives $X^2a_t=0$ for
$t=1,2$. A bracket of two basis
elements of index greater than $h$ has output index at least $2h+3$,
where $d$ vanishes, and its two action terms vanish because those
elements kill $F_1$ by \eqref{bnw:eq:degree}. Thus $d$ is an ordinary
cocycle. For every $x\in\mathfrak{w}_h$ we have $d(x)\in F_1$,
$xF_1\subseteq V$ and $xV=0$, so the restricted condition in
\eqref{bnw:eq:restricted} holds for every $x\in\mathfrak w_h$. The parameters range over a space of dimension
$2\dim F_1+(h-2)\dim V=(h+2)\dim V$, and normalised cocycles correspond
bijectively to cohomology classes.
\end{proof}

\begin{proposition}\label{bnw:prop:h2}
If $h=2$ and $p\ge11$, then $\dim H^1(u_0(\mathfrak{w}_2),L)=4\dim V$.
\end{proposition}

\begin{proof}
Put $a=d(e_3)$, $b=d(e_4)$ and $X=e_2$. The brackets
$[e_3,e_4]=e_7$ and $[e_3,e_5]=2e_8$ give
\[
\tfrac13X^2a=e_3b-e_4a,\qquad \tfrac14X^2b=e_3Xa-e_5a.
\]
Both $e_7$ and $e_8$ lie in $\mathfrak{w}_2$ because $p\ge11$. Let
$m=\max\{\operatorname{top}(a),\operatorname{top}(b)\}$, with $\operatorname{top}(0)=-\infty$, and
suppose $m\ge2$. By \eqref{bnw:eq:degree} the right side of the second
equation lies in $F_{m-3}$, while $X^2$ acts on layer $m$ with an
invertible leading map; hence $b$ has no component in layer $m$. The
right side of the first equation then also lies in $F_{m-3}$, so $a$
has no component in layer $m$ either, contradicting the choice of $m$.
Therefore $a,b\in F_1$, and the cocycles of
Proposition~\ref{bnw:prop:lower} attain the resulting bound.
\end{proof}

\begin{proposition}\label{bnw:prop:two-step}
If $h\ge4$ and $2h+5\le p\le3h+2$, then
\[
 \dim H^1(u_0(\mathfrak{w}_h),L)=(h+2)\dim V.
\]
\end{proposition}

\begin{proof}
Put $q=p-2h-2$, so that $3\le q\le h$ and $p-2=2h+q$. The ideal $\mathfrak{w}_h$ is two-step nilpotent, since the least index in
$[\mathfrak{w}_h,[\mathfrak{w}_h,\mathfrak{w}_h]]$ is $3h+1>p-2$. For $1\le j\le q$ the
brackets $[X,e_{h+j}]=j\,e_{2h+j}$ and $[X,e_{2h+j}]=0$ give
$X^2a_j=0$, so $a_j\in F_1$; for $q<j\le h$ the bracket
$[X,e_{h+j}]$ vanishes and $a_j\in V$. If $3\le j\le q$, then
$[e_{h+1},e_{h+j-1}]=(j-2)e_{2h+j}$, and the two action terms in the
corresponding cocycle equation vanish because $e_{h+1}$ and
$e_{h+j-1}$ kill $F_1$. Hence $d(e_{2h+j})=j^{-1}Xa_j=0$ and
$a_j\in V$. Only $a_1,a_2\in F_1$ and $a_3,\ldots,a_h\in V$ survive, and the
cocycles of Proposition~\ref{bnw:prop:lower} attain this bound.
\end{proof}

\begin{proposition}\label{bnw:prop:deep}
If $h\ge4$ and $p\ge3h+5$, then
$\dim H^1(u_0(\mathfrak{w}_h),L)=(h+2)\dim V$.
\end{proposition}

\begin{proof}
For $3\le s\le h$ the brackets
$[e_{h+1},e_{h+s-1}]=(s-2)e_{2h+s}$ and $[X,e_{h+s}]=s\,e_{2h+s}$
give
\begin{equation}\label{bnw:eq:forward}
\frac{s-2}{s}Xa_s=e_{h+1}a_{s-1}-e_{h+s-1}a_1,
\end{equation}
and for $t=1,2$ the brackets $[e_{h+t},e_{2h}]=(h-t)e_{3h+t}$,
$[X,e_{h+t}]=t\,e_{2h+t}$ and $[X,e_{2h+t}]=(h+t)e_{3h+t}$ give
\begin{equation}\label{bnw:eq:return}
\frac{h-t}{t(h+t)}X^2a_t=e_{h+t}a_h-e_{2h}a_t.
\end{equation}
All the basis elements displayed lie in $\mathfrak{w}_h$ because
$3h+2\le p-2$, and every coefficient shown is nonzero in $k$.

For $a_i\ne0$ put $m_i=\operatorname{top}(a_i)$ and $\omega_i=m_i+i$, and set both to
$-\infty$ when $a_i=0$; let $H=\max(\omega_1,\omega_2)$. If $m_s\ge1$, comparing
highest layers in \eqref{bnw:eq:forward} through \eqref{bnw:eq:degree} gives
\[
m_s-1\le\max\{m_{s-1}-2,\;m_1-s\},
\qquad\text{that is}\qquad
\omega_s\le\max\{\omega_{s-1},\omega_1\},
\]
while $m_s\le0$ gives $\omega_s\le s$.
Induction on $s$ now yields
\begin{equation}\label{bnw:eq:upperweight}
\omega_h\le\max\{H,h\}.
\end{equation}

Suppose $m_1\ge2$ or $m_2\ge2$. Some $t\in\{1,2\}$ then has
$m_t\ge2$ and $\omega_t=H$. Indeed, if $m_1,m_2\ge2$ take the index
attaining $H$; if $m_1\ge2$ and $m_2\le1$ then $\omega_2\le3\le \omega_1$; and
if $m_2\ge2$ and $m_1\le1$ then $\omega_1\le2<4\le \omega_2$. The left side
of \eqref{bnw:eq:return} has highest layer $m_t-2$, and the right side
lies in $F_{\max\{m_h-t-1,\,m_t-h-1\}}$. Since $m_t-h-1<m_t-2$, the
first of these must dominate, so
\[
\omega_h=m_h+h\ge m_t+t+h-1=H+h-1.
\]
Since $H\ge3$ and $h\ge4$ give $H+h-1>\max\{H,h\}$, this contradicts
\eqref{bnw:eq:upperweight}. Hence $a_1,a_2\in F_1$. Equation
\eqref{bnw:eq:forward} then gives $Xa_3=0$, and inductively $Xa_s=0$ for
every $s\ge3$, because the elements $e_j$ with $j>h$ kill $F_1$.
So $a_3,\ldots,a_h\in V$, and the cocycles of
Proposition~\ref{bnw:prop:lower} attain this bound.
\end{proof}

\begin{proof}[Proof of Theorem~\ref{thm:bnw-main}]
By \eqref{bnw:eq:BNW} it suffices to compute $\dim H^1(u_0(\mathfrak{w}_h),L)$.
The abelian range is \eqref{bnw:eq:abelian}, and $p=2h+3$ is
\eqref{bnw:eq:first}. For $p\ge2h+5$ apply Proposition~\ref{bnw:prop:h2} when
$h=2$, Proposition~\ref{bnw:prop:two-step} when $h\ge4$ and $p\le3h+2$,
and Proposition~\ref{bnw:prop:deep} when $h\ge4$ and $p\ge3h+5$. The two
integers between these last ranges are $3h+3=3(h+1)$ and $3h+4$; the
first is composite and the second even, so neither is an odd prime.
Finally $p-h-2\le h-1$ when $p\le2h+1$, while $p-h-2=h+1$ when
$p=2h+3$ and $p-h-2\ge h+3$ when $p\ge2h+5$. The three values
$p-h-2$, $h+1$ and $h+2$ therefore agree with $\min\{p-h-2,h+2\}$ on
their respective ranges.

The commutator ideal of $\mathfrak w_h$ is $\mathfrak w_{2h+1}$
when $2h+1\le p-2$, and zero otherwise. Thus the lower bound is
$h$ in the nonabelian range and $p-h-2$ in the abelian range.
The formula in \cref{thm:bnw-main} is strictly greater than $h$
in the former range and equals $p-h-2$ in the latter.
\end{proof}

\subsection{The uniserial endpoint}

\begin{proposition}[Loewy layers and higher extensions at height $p-3$]
\label{prop:w11-uniserial-higher}
Suppose $\mathfrak w=W(1;\mathbf1)$ and $h(\chi)=p-3$.
Put $N=p^{(p-1)/2}$, and let $L$ be the unique simple module.
Then
\[
 u_\chi(\mathfrak w)\cong\Mat_N(k[t]/(t^p)).
\]
The projective cover $P$ of $L$ has Loewy length $p$, with
\[
 \rad^jP/\rad^{j+1}P\cong L\qquad(0\le j<p).
\]
Its self-extension algebra is
\[
 \Ext^*_{u_\chi(\mathfrak w)}(L,L)
 \cong\Lambda_k(\alpha)\otimes_k k[\beta],
 \qquad |\alpha|=1,\quad |\beta|=2.
\]
In particular, $\dim\Ext^q(L,L)=1$ for every $q\ge0$, and
$\alpha^2=0$ although $\Ext^2(L,L)\ne0$.
\end{proposition}

\begin{proof}
By \cref{ex:w11}, $C_\chi=k[t]/(t^p)$ and the matrix factor has size $N$.
Write $C=k[t]/(t^p)$ and $J=(t)$. Then
$J^j=(t^j)$, so $J^j/J^{j+1}$ is one-dimensional for $0\le j<p$.
Morita equivalence sends the regular $C$-module to $P$, its simple
quotient to $L$, and its radical filtration to that of $P$.

A minimal projective resolution of $k=C/J$ is
\[
 \cdots\xrightarrow{t}C\xrightarrow{t^{p-1}}C
 \xrightarrow{t}C\longrightarrow k\longrightarrow0.
\]
Indeed, $\ker(t:C\to C)=(t^{p-1})$ and
$\ker(t^{p-1}:C\to C)=(t)$. Applying $\Hom_C(-,k)$ makes all
maps zero, giving a one-dimensional extension group in each degree.
For the Yoneda products, let $P_i=C$ and let $d_i$ be multiplication
by $t$ for odd $i$ and by $t^{p-1}$ for even $i$.
The degree-two lift $P_{i+2}\to P_i$ given by the identity represents
$\beta$; its powers implement the two-periodicity. A degree-one
comparison lift representing $\alpha$ has components
\[
 A_i:P_{i+1}\to P_i,\qquad
 A_i=\begin{cases}1,&i\text{ even},\\t^{p-2},&i\text{ odd}.
 \end{cases}
\]
They satisfy $d_iA_i=A_{i-1}d_{i+1}$ for $i\ge1$.
The augmentation composed with $A_0A_1$ is zero, because $p>2$,
so $\alpha^2=0$. The degree-two lift commutes with these components,
which repeat after two steps. Thus $\alpha\beta=\beta\alpha$,
and $\beta^j$ and $\alpha\beta^j$ give the nonzero classes in every
even and odd degree, respectively. Morita equivalence preserves
Yoneda products and gives the stated algebra.
\end{proof}

\subsection{A primitive corner at height two}

\begin{proposition}[A corner at height two]
\label{prop:w11-height-two-corner}
Let $\mathfrak w=W(1;\mathbf1)$ and normalize a character of height two
by $\chi(e_1)=1$ and $\chi(e_i)=0$ for $i\ne1$. Put
$A=u_\chi(\mathfrak w)$, $F=e_{-1}$, $H=e_0$, $E=e_1$, and
\[
\begin{aligned}
 e&=1-H^{p-1},&\qquad Y&=eEF e,\\
 T&=Y-\tfrac14e,& q&=\left(T^{(p-1)/2}-(-1)^{(p-1)/2}e\right)^2.
\end{aligned}
\]
Then $q$ is a primitive idempotent of $A$, and $Aq$ is the projective
cover of its unique simple module $L$. Moreover,
\[
 A\cong\Mat_{p^2}(qAq),\qquad \dim qAq=p^{p-4}.
\]
For $\mathfrak l=\langle F,H,E\rangle\cong\mathfrak{sl}_2$,
multiplication identifies
\[
 Aq\cong A\otimes_{u_\chi(\mathfrak l)}u_\chi(\mathfrak l)q,
\]
where $u_\chi(\mathfrak l)q$ is a simple projective module of dimension $p$.
\end{proposition}

\begin{proof}
In $u_\chi(\mathfrak l)$ one has $E^p=1$, $F^p=0$ and $H^p=H$.
The element $K=H^2-H-EF$ commutes with $E,F,H$. Since
$Ef(H)=f(H-1)E$, induction gives
\[
 0=E^pF^p=\prod_{a\in\mathbf F_p}
 \bigl((H-a)^2-(H-a)-K\bigr).
\]
Multiplying by $e$ and using $eH=0$ and $eK=-Y$ yields
\[
 \prod_{a\in\mathbf F_p}\bigl(Y+(a^2+a)e\bigr)=0.
\]
With $T=Y-\tfrac14e$ and $m=(p-1)/2$, this reads
\[
 T\bigl(T^m-(-1)^me\bigr)^2=0.
\]
The root $T=0$ is simple, and the remaining roots are double.
The polynomial defining $q$ takes value one at zero and vanishes to
order two at every other root, so $q^2=q$.
The matrix units for $H,E$ identify $u_\chi(\mathfrak l)$ with
$\Mat_p(eu_\chi(\mathfrak l)e)$. Its corner has basis
$e,Y,\ldots,Y^{p-1}$ by the PBW theorem. Thus $q$ selects its simple
root $Y=\tfrac14e$, and $u_\chi(\mathfrak l)q$ is simple and projective
of dimension $p$.

Use the induced realization of $L$ in \cref{subsec:w11-even} at $h=2$,
with generator $v$. Its $H$-weight-zero subspace has basis
\[
 w_a=F^a\bigl(1-(H-a)^{p-1}\bigr)v\qquad(0\le a<p).
\]
The brackets in $\mathfrak l$ give
\[
 Yw_a=w_{a+1}-a(a+1)w_a,\qquad w_p=0.
\]
The eigenvalue $\tfrac14$ occurs once, at $a=(p-1)/2$.
Consequently $\dim qL=1$. Since $A$ has a unique simple module of
dimension $p^2$, the idempotent $q$ is primitive, $Aq$ is its projective
cover, and $A\cong\Mat_{p^2}(qAq)$. The corner dimension follows from
$\dim A=p^p$. Finally $q\in u_\chi(\mathfrak l)$, so multiplication gives
the displayed induction isomorphism.
\end{proof}

\subsection{A central element of the height-two corner}

\begin{proposition}\label{prop:w11-central-transfer}
Put $f=e_{-1}$ and $m=(p+1)/2$. The element
\[
 Z_p=(\ad f)^{p-1}(e_{p-2}^{\,m})\in U(W(1;\mathbf1))
\]
is central and has nonzero image in every reduced enveloping algebra.
At height two, with $q$ as in \cref{prop:w11-height-two-corner},
\[
 0\ne qZ_pq\in Z(qAq)\cap\Jac(qAq).
\]
\end{proposition}

\begin{proof}
Write $v=e_{p-2}^{\,m}$ and $D=\ad f$.
For the adjoint action, $e_iv=0$ for $i\ge1$ and $e_0v=-v$,
since $m(p-2)=-1$ in $k$. The centrality of $f^p$ gives $D^p=0$.
Thus $fZ_p=e_0Z_p=0$ under the adjoint action.
For $1\le i\le p-2$, commuting $e_i$ through $D^{p-1}$ leaves
only the terms with $e_0$ and $f$ acting on $v$:
\[
 e_iD^{p-1}v
 =(-1)^{i+1}(i+1)!
 \left(\binom{p-1}{i}+\binom{p-1}{i+1}\right)D^{p-1-i}v=0.
\]
Here $\binom{p-1}{j}=(-1)^j$ in characteristic $p$.

In the PBW symbol, the monomial
$e_{p-4}^{m-2}e_{p-3}^2$ has coefficient
\[
 \binom m2\frac{(p-1)!}{(2!)^{m-2}}
 \bigl((p-1)(p-2)\bigr)^{m-2}(p-1)^2\ne0.
\]
It comes from differentiating two factors once and the other $m-2$
factors twice. The PBW degree is $m<p$, so this symbol survives
unchanged in every reduced enveloping algebra.

At height two, \eqref{bnw:eq:degree} gives
$e_{p-2}F_a\subseteq F_{a-(p-3)}$.
Since $m(p-3)>p-1$, the operator $v$ annihilates $L$;
its iterated commutator $Z_p$ therefore annihilates $L$ as well.
Thus $Z_p\in\Jac(A)$, since $L$ is the unique simple module.
The idempotent $q$ is full and $Z_p$ is central, so $qZ_pq=0$
would imply $Z_p=0$ by the matrix reconstruction formula.
This proves the last assertion.
\end{proof}

\subsection{Equivariant resolutions and invariants}

\begin{lemma}\label{lem:tate-invariants}
Let $D=k[X]/(X^p)$, let $R=B\rtimes D$ be a restricted enveloping
algebra with $B=u_0(\mathfrak m)$, and let $\mathcal T$ be an
$R$-equivariant resolution of $k$ by free $B$-modules.
For an $R$-module $L$, put $K=\Hom_B(\mathcal T,L)$.
If every $K^q$ is free over $D$, then
\[
 H^*(R,L)\cong H^*(K^D).
\]
\end{lemma}

\begin{proof}
Let $\mathcal P$ be the periodic free $D$-resolution of $k$,
with maps alternating between $X$ and $X^{p-1}$.
Inflate it to $R$ by the augmentation of $B$.
The diagonal Hopf action makes $\mathcal T\otimes_k\mathcal P$
a free $R$-resolution. To see freeness, for a $B$-free $R$-module $M$
use the triangular isomorphism
\[
 R\otimes_B M\longrightarrow M\otimes_k D,\qquad
 bX^r\otimes v\longmapsto
 \sum_{j=0}^r\binom rj bX^jv\otimes X^{r-j}.
\]
Exactness follows by taking homology in the two tensor directions.
Its cochains form $\Hom_D(\mathcal P,K)$.
Since $D$ is Frobenius, its free modules are injective;
taking $D$-cohomology first leaves $K^D$ in column zero.
\end{proof}

\subsection{An equivariant resolution in characteristic seven}

\begin{theorem}[Higher self-extensions in characteristic seven]
\label{thm:w11-height-two-all-ext}
Suppose $p=7$ and $\chi$ has height two. For the unique simple
$A=u_\chi(W(1;\mathbf1))$-module $L$ and every $n\ge0$,
\[
 \dim\Ext_A^n(L,L)=\binom{n+2}{2},\qquad
 \sum_{n\ge0}\dim\Ext_A^n(L,L)t^n=\frac1{(1-t)^3}.
\]
The minimal projective resolution of the simple module of the corner
in \cref{prop:w11-height-two-corner} has rank $\binom{n+2}{2}$ in degree $n$.
\end{theorem}

\begin{proof}
Normalize $\chi$ as in \cref{prop:w11-height-two-corner}. Put
$\mathfrak n=\langle e_2,e_3,e_4,e_5\rangle$,
$\mathfrak m=\langle e_3,e_4,e_5\rangle$, and $X=e_2$.
The only nonzero bracket in this basis is $[X,e_3]=e_5$;
all restricted seventh powers vanish. Consequently
\[
 R=u_0(\mathfrak n)=B\rtimes D,\quad
 B=k[b_3,b_4,b_5]/(b_3^7,b_4^7,b_5^7),\quad D=k[X]/(X^7),
\]
where $b_i=e_i$ and $X$ acts on $B$ by the derivation
$\delta(b_3)=b_5$, $\delta(b_4)=\delta(b_5)=0$.
Lemma~\ref{bnw:lem:pivot} gives that $L$ is free of rank seven over
$D$, with $\ker X^2=F_1$. Equation~\eqref{bnw:eq:degree} gives
$\mathfrak m F_1=0$, $b_3^5L=0$, and $b_i^6L=0$ for $3\le i\le5$.

Use the following free $B$-resolution of $k$:
\[
 \mathcal T=B\otimes\Lambda(\epsilon_3,\epsilon_4,\epsilon_5)
 \otimes\Gamma(\rho_3,\rho_4,\rho_5),\qquad
 d\epsilon_i=b_i,\quad d\rho_i=b_i^6\epsilon_i.
\]
Here $\deg\epsilon_i=1$, $\deg\rho_i=2$, and $\Gamma$ denotes the
 divided-power algebra. This is the tensor product of the three
one-variable resolutions whose maps alternate between $b_i$ and $b_i^6$,
so it resolves $k$. Extend $\delta$ to it by
\[
 \delta\epsilon_3=\epsilon_5,\quad
 \delta\rho_3=b_3^5\epsilon_3\epsilon_5,\quad
 \delta\epsilon_4=\delta\epsilon_5=\delta\rho_4=\delta\rho_5=0,
\]
using the derivation rule, including
$\delta\gamma_r(\rho_i)=(\delta\rho_i)\gamma_{r-1}(\rho_i)$.
The identity $d\delta=\delta d$ on $\rho_3$ reads
\[
 d(b_3^5\epsilon_3\epsilon_5)
 =b_3^6\epsilon_5-b_3^5b_5\epsilon_3
 =\delta(b_3^6\epsilon_3).
\]
It holds on the other generators as well. Also $\delta^7=0$.
For divided powers of $\rho_3$, every iterated derivative containing
 two factors $\delta^j\rho_3$ vanishes because each such factor contains
$\epsilon_3\epsilon_5$. The remaining seventh derivative is zero
since it differentiates $b_3^5$ six times. Thus $\mathcal T$ is
an $R$-equivariant $B$-resolution.

Put $K=\Hom_B(\mathcal T,L)$, with
$(Xf)(v)=Xf(v)-f(\delta v)$. The terms $b_i^6L=0$ remove the
 divided-power part of the cochain differential, and $b_3^5L=0$
removes the correction $\delta\rho_3$ from the $D$-action. Hence
\[
 K^r=\bigoplus_{a+2b=r}
 L\otimes\Lambda^a(\mathfrak m^*)\otimes S^b(k^3)
\]
as $D$-modules; $D$ acts trivially on the last factor and by the
coadjoint action on the exterior factor. A diagonal tensor product
of a free $D$-module with any $D$-module is free, by the same triangular
Hopf map. Each $K^r$ is therefore free. By \cref{lem:tate-invariants},
$H^*(R,L)=H^*(K^D)$.

The coadjoint operator $N$ on $\Lambda^a(\mathfrak m^*)$ has square
zero: it replaces the dual of $e_5$ by the dual of $e_3$ and kills
 the other basis vectors. If a cochain is $D$-invariant, its coefficient
vector satisfies $Xf=-Nf$, so $X^2f=0$. All coefficients belong to
$F_1$, which $\mathfrak m$ kills. The remaining Koszul differential
therefore vanishes on $K^D$. Each free $D$-summand has one-dimensional
invariants per copy of $D$, giving
\[
 \sum_{r\ge0}\dim H^r(R,L)t^r
 =7\frac{(1+t)^3}{(1-t^2)^3}=\frac7{(1-t)^3}.
\]

Since $Q\cong\Mat_7(k)$, \cref{bnw:prop:adjunction} gives
\[
 \dim\Ext_A^n(L,L)=\tfrac17\dim H^n(R,L)=\binom{n+2}{2}.
\]
For the local basic corner, these dimensions are the ranks of the
minimal free resolution.
\end{proof}

\subsection{Higher extensions in the abelian range}

\begin{proposition}[Higher self-extensions in the abelian range]
\label{prop:w11-abelian-all-ext}
Let $h$ be even, $2\le h\le p-3$, and $p\le2h+1$.
Put $e=p-h-2$. For the unique simple
$A=u_\chi(W(1;\mathbf1))$-module $L$ at height $h$,
\[
 \sum_{n\ge0}\dim\Ext_A^n(L,L)t^n=\frac1{(1-t)^e},\qquad
 \dim\Ext_A^n(L,L)=\binom{n+e-1}{e-1}.
\]
\end{proposition}

\begin{proof}
Normalize $\chi$ as in \cref{subsec:w11-even} and put $h=2s$.
The inequality $p\le2h+1$ makes $\mathfrak w_h$ abelian:
the smallest bracket of distinct basis vectors has index $2h+1>p-2$.
Write $X=e_h$, $\mathfrak m=\langle e_{h+1},\ldots,e_{p-2}\rangle$,
$D=k[X]/(X^p)$ and $B=u_0(\mathfrak m)$.
Thus $u_0(\mathfrak w_h)=B\otimes D$, and $B$ is a truncated
polynomial algebra on $e$ generators. By Lemma~\ref{bnw:lem:pivot},
$L$ is free of rank $p^s$ over $D$ and $L^D=\ker X=V$.
Equation~\eqref{bnw:eq:degree} gives $\mathfrak m V=0$ and
$e_i^{p-1}L=0$ for every $i\ge h+1$, since each such action lowers
the filtration by at least two.

Take the tensor product Tate resolution $\mathcal T$ of $k$ over $B$,
with $e$ exterior generators in degree one and $e$ divided-power
generators in degree two. Give these generators trivial $D$-action.
In $K=\Hom_B(\mathcal T,L)$ the
 divided-power part of the differential is zero, since its coefficients
are the operators $e_i^{p-1}$.
Every $K^r$ is a direct sum of copies of the free $D$-module $L$,
so \cref{lem:tate-invariants} identifies cohomology with that of $K^D$. Its coefficients belong to
$V$, annihilated by $\mathfrak m$, so its remaining differential is zero.
Counting exterior and divided-power monomials gives
\[
 \sum_{r\ge0}\dim H^r(u_0(\mathfrak w_h),L)t^r
 =p^s\frac{(1+t)^e}{(1-t^2)^e}=\frac{p^s}{(1-t)^e}.
\]
Apply \cref{bnw:prop:adjunction} to divide by $p^s$ in every degree.
\end{proof}

\subsection{Central operators and the Yoneda algebra}

\begin{proposition}[Central polynomial operators]
\label{prop:w11-p7-central-operators}
Under the hypotheses of \cref{thm:w11-height-two-all-ext}, there are
central classes $\theta_3,\theta_4,\theta_5$ of degree two such that
\[
 P=k[\theta_3,\theta_4,\theta_5]\subseteq
 E^*=\Ext_A^*(L,L)
\]
is a polynomial algebra. As a graded $P$-module, $E^*$ is free,
with one basis element in degree zero, three in degree one,
three in degree two, and one in degree three.
\end{proposition}

\begin{proof}
Put $\mathfrak g=W(1;\mathbf1)$ and let $\lambda_i$ be the coordinate
functional dual to $e_i$, for $i=3,4,5$. The split ordinary Lie
extension $\mathfrak g\oplus kc$, with $c$ central, becomes a
restricted extension on setting
\[
 c^{[7]}=0,\qquad
 (x+ac)^{[7]}=x^{[7]}+\lambda_i(x)^7c.
\]
The added term is additive and seven-semilinear and has zero adjoint
action, so all restricted identities hold. Its restricted
cohomology class is $\gamma_i\in H^2(u_0(\mathfrak g),k)$.

Tensoring an extension representing $\gamma_i$ with an $A$-module
defines a natural degree-two extension operator: the coproduct
identifies the restricted-power relation on the tensor product with
the sum of the relations on its factors, so tensoring a zero-character
module with a character-$\chi$ module retains character $\chi$.
The symmetry of tensor products and the interchange law for
composition of extensions make this operator commute with every
Yoneda class \cite[Section 9.3]{Witherspoon2019}.
The same interchange argument applies to the category of
character-$\chi$ modules.
Let $\theta_i$ be its value on $L$.

We identify its action on the resolution in
\cref{thm:w11-height-two-all-ext}. On the divided-power factor define
$\partial_{\rho_i}$ by
\[
 \partial_{\rho_i}\gamma_r(\rho_i)=\gamma_{r-1}(\rho_i)
 \quad(r\ge1),\qquad \partial_{\rho_i}(1)=0,
\]
acting trivially on the other factors. It has homological degree
minus two. It commutes with the differential, since every
$d\rho_j=b_j^6\epsilon_j$ is independent of the $\rho$ variables.
It also commutes with the $X$-action, since its sole correction
$\delta\rho_3=b_3^5\epsilon_3\epsilon_5$ is independent of them.
Thus $\partial_{\rho_i}\otimes1$ is an $R$-linear comparison map
on the total free resolution.

The augmentation of this comparison map has zero ordinary
exterior values. Its restricted-power value is one on $e_i$ and zero
on the other basis vectors. It therefore represents the restriction
of $\gamma_i$ to $R=u_0(\mathfrak w_2)$, and its cochain action is
multiplication by the polynomial variable dual to $\rho_i$.
Induction adjunction respects this action, by the tensor identity
\[
 U\otimes\operatorname{Ind} S
 \cong\operatorname{Ind}(\operatorname{Res}U\otimes S).
\]
The semisimple quotient reduction is natural as well. Thus the
identification
\[
 E^*=\Hom_Q(S,H^*(R,L))
\]
is compatible with all three operators. In particular their actions
on $H^*(R,L)$ are $Q$-linear.

Let $P$ be the polynomial ring on the three comparison operators.
The invariant-complex proof of
\cref{thm:w11-height-two-all-ext} identifies $H^*(R,L)$, as a
$P$-module, with a free module whose basis has dimensions
$7,21,21,7$ in degrees zero through three.
As $Q=\Mat_7(k)$, a primitive
matrix idempotent identifies $\Hom_Q(S,H^*(R,L))$ with a direct
summand of this free $P$-module. Hence $E^*$ is a finite graded
projective $P$-module, and the sequence
$\theta_3,\theta_4,\theta_5$ is regular on it.
Taking successive quotients therefore gives
\[
 \sum_{n\ge0}\dim(E^*/P_+E^*)_nt^n
 =(1-t^2)^3\sum_{n\ge0}\dim E^nt^n=(1+t)^3.
\]
Choose homogeneous lifts of a basis of this quotient.
Graded Nakayama gives a surjection from the free $P$-module on these
eight lifts to $E^*$. Its source and target have the same Hilbert
series by \cref{thm:w11-height-two-all-ext}, so it is an isomorphism
in every degree. In particular the action of $P$ on its nonzero
unit class is faithful. This proves that the $\theta_i$ generate
a polynomial subalgebra of $E^*$ and proves the asserted freeness.
\end{proof}

\begin{theorem}[The Yoneda algebra at height two in characteristic seven]
\label{thm:w11-p7-yoneda}
Under the hypotheses of \cref{thm:w11-height-two-all-ext}, let
\[
 B=k\langle a,b,c\rangle/
 (a^2,\ ba+ab,\ b^2-5ab,\ ca-3ab+ac,\
 cb-2ab-4ac+bc,\ c^2-5bc),
\]
with $a,b,c$ in degree one. Then
\[
 \Ext_A^*(L,L)\cong B\otimes k[\theta_3,\theta_4,\theta_5],
 \qquad \deg\theta_i=2,
\]
where the $\theta_i$ are the central operators of
\cref{prop:w11-p7-central-operators}.
The algebra $B$ is symmetric of dimension eight, with basis
$1,a,b,c,ab,ac,bc,abc$ and Hilbert series $(1+t)^3$.
The square-zero elements of $B^1$ form the line $ka$.
\end{theorem}

\begin{proof}
Write
\[
 w_{r,j}=F^r\bigl(1-(H-j)^6\bigr)v\qquad
 (0\le r<7,\ j\in\mathbf F_7).
\]
Subscripts $j$ are read modulo seven; $w_{r,j}=0$ for $r<0$ or $r>6$.
The induced realization of $L$ in \cref{subsec:w11-even} gives
\[
\begin{aligned}
 Fw_{r,j}&=w_{r+1,j},&
 Hw_{r,j}&=(j-r)w_{r,j},\\
 Ew_{r,j}&=w_{r,j+1}-r(2j-r+1)w_{r-1,j},\\
 Xw_{r,j}&=-3r\,w_{r-1,j+1}
 +\left(6j\binom r2-6\binom r3\right)w_{r-2,j}.
\end{aligned}
\]
Here $X=e_2$; the other actions are obtained by
$e_i=(i-2)^{-1}[E,e_{i-1}]$ for $i=3,4,5$.

Define endomorphism-valued cocycles $D_a,D_b,D_c$ to vanish on
$F,H,E$. Their values on $X$ are
\[
\begin{aligned}
 D_a(X)w_{r,j}&=(5-j)w_{r,j+2},\\
 D_b(X)w_{r,j}&=(2j^2-j+4)w_{r,j+2},\\
 D_c(X)w_{r,j}&=(4j^3+4j)w_{r,j+2}-j\,w_{r+1,j+3};
\end{aligned}
\]
their values on $e_3,e_4,e_5$ are defined by the same bracket recursion
with $E$. Substitution in the preceding action formulas gives
\[
 D([e_i,e_j])=[e_i,D(e_j)]+[D(e_i),e_j],\qquad
 \sum_{q=0}^6 e_i^qD(e_i)e_i^{6-q}=D(e_i^{[7]}).
\]
Thus these maps represent restricted self-extensions. Normalizing
their restrictions to $\mathfrak n=\langle e_2,e_3,e_4,e_5\rangle$
to have value zero on $X$, their values on $e_3,e_4$ have parameters
\[
 (1,0,0),\quad(6,0,2),\quad(6,1,4)
\]
in the form
$d(e_3)=\alpha w_{0,3}+\beta w_{1,4}$,
$d(e_4)=\gamma w_{0,4}$; also $d(e_5)=Xd(e_3)$.
These three triples are independent. By
\cref{thm:bnw-main} the associated classes $a,b,c$ form a basis of $E^1$.

Use restricted degree-two cochains $(\phi,\omega)$ on
$\mathfrak n$, with $\phi$ alternating and $\omega$ specified on its
four basis vectors. Their coboundaries are
\[
 (\delta f)(x,y)=x f(y)-y f(x)-f([x,y]),\qquad
 \omega_f(x)=x^6f(x).
\]
The cup of $D_\alpha,D_\beta$, evaluated on $w_{0,0}$ under
induction adjunction, has
\[
 \phi_{\alpha\beta}(x,y)=
 \bigl(D_\alpha(x)D_\beta(y)-D_\alpha(y)D_\beta(x)\bigr)w_{0,0}.
\]
Its power component is a sum of terms
\[
 x^uD_\alpha(x)x^vD_\beta(x)x^w w_{0,0},\qquad u+v+w=5.
\]
Each cocycle raises the $F$-filtration by at most one, whereas
$x\in\mathfrak n$ lowers it by at least one, so every term is zero.

The following table supplies coboundaries for the six stated
relations. In the table $t_{x,r}=w_{r,x+r}$, and all unlisted values,
including $f(e_5)$, are zero.
\[
\begin{array}{c|lll}
 \text{relation}&f(e_2)&f(e_3)&f(e_4)\\ \hline
 a^2&2t_{2,2}&0&0\\
 ba+ab&4t_{2,2}+5t_{2,3}&5t_{3,2}&0\\
 b^2-5ab&t_{2,2}+t_{2,3}&t_{3,2}&0\\
 ca-3ab+ac&6t_{2,2}+4t_{2,3}+2t_{2,4}&5t_{3,3}&3t_{4,2}\\
 cb-2ab-4ac+bc&5t_{2,2}+2t_{2,3}+3t_{2,4}
   &2t_{3,2}+4t_{3,3}&0\\
 c^2-5bc&t_{2,3}+2t_{2,4}+6t_{2,5}
   &t_{3,2}+3t_{3,3}+3t_{3,4}&t_{4,2}+t_{4,3}
\end{array}
\]
The identities $\phi_{\mathrm{relation}}=\delta f$ follow directly
from the four action formulas above. Every displayed $f(e_i)$ lies
in $F_5$, so all $\omega_f(e_i)$ vanish. These are therefore
restricted coboundaries.

For a weight-zero cochain, let $\phi_{xy,r}$ be the coefficient
of $w_{r,x+y+r}$ in $\phi(e_x,e_y)$, and let $\omega_{x,r}$ be
the coefficient of $w_{r,r}$ in $\omega(e_x)$.
The following table gives the coefficients of three linear forms
$\ell_1,\ell_2,\ell_3$; omitted entries are zero.
\[
\begin{array}{c|rrr@{\qquad}c|rrr}
 &\ell_1&\ell_2&\ell_3&&\ell_1&\ell_2&\ell_3\\ \hline
\phi_{23,2}&1&0&0&\phi_{23,3}&1&5&4\\
\phi_{23,4}&5&6&2&\phi_{23,5}&5&6&2\\
\phi_{23,6}&4&2&3&\phi_{24,1}&1&0&0\\
\phi_{24,2}&1&6&5&\phi_{24,3}&2&4&1\\
\phi_{24,4}&1&2&4&\phi_{24,5}&6&5&3\\
\phi_{25,1}&4&4&0&\phi_{25,2}&2&1&5\\
\phi_{25,4}&2&1&5&\phi_{34,0}&4&0&0\\
\phi_{34,1}&1&5&3&\phi_{35,0}&5&2&0\\
\omega_{2,0}&3&1&0&&&&\\
\end{array}
\]
Substitution in the coboundary formula shows that these forms
annihilate every weight-zero coboundary. Their values on $ab,ac,bc$
are the three coordinate vectors. They also vanish on
$\theta_3,\theta_4,\theta_5$, whose cochains have
$\phi=0$ and sole power value $\omega(e_i)=w_{0,0}$.
Toral weight projection shows that a boundary for a weight-zero
cochain can always be chosen with weight zero.
Consequently
\[
 E^2=\langle ab,ac,bc\rangle\oplus
 \langle\theta_3,\theta_4,\theta_5\rangle.
\]

It remains to test the top product. In degree three write cochains
as $(\psi,\eta)$, with ordinary values on the four increasing
triples and power-compatibility values $\eta(x,y)$ for ordered
$x,y\in\{e_2,e_3,e_4,e_5\}$. For degree-two boundaries the latter are
\begin{equation}\label{eq:p7-power-compatibility}
 R_{\phi,\omega}(x,y)=
 x^6\phi(x,y)+x^5\phi(x,[x,y])+y\omega(x).
\end{equation}
Here $(\operatorname{ad}x)^2\mathfrak n=0$.
These boundary tests can also be read directly from a free comparison
complex with generators $\epsilon_x$ in degree one,
$\epsilon_{xy},\rho_x$ in degree two, and
$\epsilon_{xyz},\rho_x\epsilon_y$ in degree three. Its first maps are
\[
\begin{aligned}
 d\epsilon_x&=x,\\
 d\epsilon_{xy}&=x\epsilon_y-y\epsilon_x-\epsilon_{[x,y]},\\
 d\rho_x&=x^6\epsilon_x,\\
 d(\rho_x\epsilon_y)&=
 x^6\epsilon_{xy}+x^5\epsilon_{x,[x,y]}+y\rho_x;
\end{aligned}
\]
the map on $\epsilon_{xyz}$ is the ordinary Jacobi map.
In the normalized bar resolution its comparison maps send
$\epsilon_{xy}$ to $[x|y]-[y|x]$, $\rho_x$ to $[x^6|x]$,
and $\epsilon_{xyz}$ to the alternating bar sum.
Writing $v=[x,y]$, the degree-three power generator is sent to
\[
\begin{aligned}
 &[x^6|x|y]-[x^6|y|x]+[x^5|x|v]\\
 &\hspace{25mm}-[x^5|v|x]+[y|x^6|x].
\end{aligned}
\]
The bar differential gives the displayed comparison identities;
a term with a zero bar entry is interpreted as zero.
Thus any coboundary pulls back to the ordinary differential of
$\phi$ and the values in \eqref{eq:p7-power-compatibility}.
This suffices for the boundary tests below.

The cup $abc$ has power components zero and ordinary values
\[
 (\psi_{234},\psi_{235},\psi_{245},\psi_{345})
 =(6w_{1,3},\,6w_{0,3},\,3w_{0,4},\,w_{0,5}).
\]
These are the alternating sums of
$D_a(x)D_b(y)D_c(z)w_{0,0}$.
On each displayed power-comparison term, the cup expands into
products with three cocycle factors and at least four action factors.
Their filtration bound is at most $3-4<0$, giving zero.

For a weight-zero degree-three cochain, write $\psi_{xyz,r}$
for the coefficient of $w_{r,x+y+z+r}$ in $\psi(e_x,e_y,e_z)$,
and $\eta_{xy,r}$ for the coefficient of $w_{r,y+r}$ in
$\eta(e_x,e_y)$. The linear form
\[
\begin{aligned}
 \ell={}&2\psi_{234,3}+2\psi_{234,4}
 +2\psi_{235,2}+4\psi_{235,3}+3\psi_{235,4}+\psi_{235,5}\\
 &+2\psi_{245,1}+4\psi_{245,2}+\psi_{345,0}\\
 &+6\eta_{22,2}+6\eta_{22,3}+5\eta_{22,4}
 +3\eta_{23,1}+3\eta_{24,0}
\end{aligned}
\]
annihilates every degree-two boundary, by substitution in
\eqref{eq:p7-power-compatibility} and the ordinary differential.
It takes value one on $abc$.
Its power coefficients use only pairs with first input $e_2$.
On the other hand, a product $\theta_i d$ has zero ordinary part
and power values $\eta(x,y)=\delta_{x,e_i}d(y)$.
Thus $\ell$ vanishes on
$\langle\theta_3,\theta_4,\theta_5\rangle E^1$.
We obtain
\[
 E^3=k\,abc\oplus
 \langle\theta_3,\theta_4,\theta_5\rangle E^1,
\]
since \cref{prop:w11-p7-central-operators} makes the second
summand nine-dimensional.

Orient the six quadratic relations with leading words
$aa$, $ba$, $bb$, $ca$, $cb$ and $cc$, using degree-lexicographic order
with $c>b>a$.
The ten overlaps, in the order
\[
 aaa,\ baa,\ bba,\ bbb,\ caa,\ cba,\ cbb,\ cca,\ ccb,\ ccc,
\]
reduce in either order to
\[
 0,0,0,0,0,6abc,5abc,5abc,4abc,6abc.
\]
The resulting normal words are
$1,a,b,c,ab,ac,bc,abc$, so they form a basis of $B$.
The preceding cochain tests show that their images form a basis
of $E^*/P_+E^*$.
By \cref{prop:w11-p7-central-operators} their lifts are a
$P$-basis of $E^*$. Centrality of $P$ now makes multiplication
$B\otimes P\to E^*$ an algebra isomorphism in every degree.

Finally $B^2$ is central and all words of length four vanish.
The coefficient of $abc$ defines a symmetric nondegenerate pairing:
between $(a,b,c)$ and $(ab,ac,bc)$ its matrix is
\[
 \begin{pmatrix}0&0&1\\0&6&5\\1&5&4\end{pmatrix},
\]
with determinant one; the degree-zero/degree-three pairing is
nondegenerate as well. The square of $xa+yb+zc$ is
\[
 (5y^2+3xz+2yz)ab+4yz\,ac+5z^2bc.
\]
It vanishes exactly when $y=z=0$, giving the asserted line.
\end{proof}

\begin{corollary}\label{cor:w11-p7-yoneda-center}
Let $E^*$ be the algebra of \cref{thm:w11-p7-yoneda} and
$P=k[\theta_3,\theta_4,\theta_5]$. Its center and graded center coincide:
\[
 Z(E^*)=Z_{\mathrm{gr}}(E^*)
 =P\oplus P\,ab\oplus P\,ac\oplus P\,bc\oplus P\,abc.
\]
The ideal spanned by the last four summands has square zero.
Every even-degree class is central. The ideal
$N=B_{>0}\otimes P$ satisfies $N^4=0$ and $E^*/N=P$;
it consists exactly of the nilpotent elements of $E^*$.
The $P$-linear functional extracting the coefficient of $abc$ is
a symmetrizing form on the free $P$-algebra $E^*$.
\end{corollary}

\begin{proof}
The relations give $Z(B)=k\oplus B^2\oplus B^3$.
For the only possible missing component, commuting a linear
combination of $a,b,c$ with $b$ forces its $a,c$ coefficients to
vanish, and commuting the remainder with $a$ forces its $b$
coefficient to vanish. For the graded center, a degree-one
combination must first have square zero, hence belong to $ka$;
the identity $ac+ca=3ab\ne0$ then excludes it.
Since $P$ is central, comparison of polynomial coefficients
gives both center formulas. Products of their positive $B$-degree
summands have degree at least four and hence vanish.

The ideal $N$ has fourth power zero, and any element outside $N$
has nonzero image in the polynomial domain $P$ and cannot be
nilpotent. The last assertion follows by extending the symmetric
nondegenerate pairing on $B$ from \cref{thm:w11-p7-yoneda}
over $P$.
\end{proof}

\section*{Acknowledgements}

I thank J\"org Feldvoss for suggesting that the general results be tested
on small cases; the examples for $W(2)$, $H(2)^{(2)}$ and
$W(1;1)$ grew out of that suggestion.


\begin{thebibliography}{99}

\bibitem{BenkartFeldvoss2015}
G.~Benkart and J.~Feldvoss,
\emph{Some problems in the representation theory of simple modular Lie
algebras},
in: Lie Algebras and Related Topics,
Contemp. Math. 652, American Mathematical Society, Providence, RI, 2015,
207--228.

\bibitem{BoeNakanoWiesner2009}
B.~D.~Boe, D.~K.~Nakano, and E.~Wiesner,
\emph{$\Ext^1$-quivers for the Witt algebra $W(1,1)$},
J. Algebra 322 (2009), no.~5, 1548--1564.

\bibitem{Chang1941}
H.-J.~Chang,
\emph{\"Uber Wittsche Lie-Ringe},
Abh. Math. Sem. Univ. Hamburg 14 (1941), 151--184.

\bibitem{FeldvossNakano1998}
J.~Feldvoss and D.~K.~Nakano,
\emph{Representation theory of the Witt algebra},
J. Algebra 203 (1998), no.~2, 447--469.

\bibitem{Ma2026Divergence}
C.~Ma,
\emph{The Derived Algebra of Nonlinear Polynomial Divergence-Free Vector
Fields}, arXiv:2610.01838 (2026).

\bibitem{Mygind2014}
M.~Mygind,
\emph{Restricted Cartan type Lie algebras and the coadjoint representation},
Transform. Groups 20 (2015), no.~2, 495--505.

\bibitem{Rian2011}
K.~Rian,
\emph{Extensions of the Witt algebra and applications},
J. Algebra Appl. 10 (2011), no.~6, 1233--1259.

\bibitem{RianThesis2011}
K.~Rian,
\emph{Representations of the Witt--Jacobson Lie algebras},
Ph.D. thesis, Aarhus University, 2011.

\bibitem{StradeFarnsteiner1988}
H.~Strade and R.~Farnsteiner,
\emph{Modular Lie Algebras and Their Representations},
Monographs and Textbooks in Pure and Applied Mathematics, vol.~116,
Marcel Dekker, New York, 1988.

\bibitem{Witherspoon2019}
S.~J.~Witherspoon,
\emph{Hochschild Cohomology for Algebras},
Graduate Studies in Mathematics \textbf{204},
American Mathematical Society, Providence, RI, 2019.

\end{thebibliography}
\end{document}